\documentclass[11pt]{amsart}
\usepackage[mathscr]{eucal}
\usepackage{mathrsfs}
\usepackage{amsfonts}
\usepackage{amsmath}
\usepackage{amsthm}
\usepackage{amssymb}
\usepackage{latexsym}
\usepackage[all]{xy}

\theoremstyle{plain}
\newtheorem{definition}[equation]{Definition}
\newtheorem{corollary}[equation]{Corollary}
\newtheorem{lemma}[equation]{Lemma}
\newtheorem{proposition}[equation]{Proposition}
\newtheorem{theorem}[equation]{Theorem}

\theoremstyle{definition}
\newtheorem{remark}[equation]{Remark}
\newtheorem{example}[equation]{Example}

\numberwithin{equation}{section}

\def\A{{\mathscr{A}}}
\def\B{{\mathscr{B}}}
\def\C{{\mathbb{C}}}
\def\CC{{\mathscr{C}}}
\def\E{{\mathscr{E}}}
\def\k{{\Bbbk}}
\def\N{{\mathscr{N}}}
\def\QQ{{\overline{Q}}}
\def\X{{\mathcal{X}}}
\def\Z{{\mathbb{Z}}}

\newcommand{\al}{{\alpha}}
\newcommand{\g}{{\gamma}}
\newcommand{\GG}{{{\mathbf{\Gamma}}_n}}
\newcommand{\hh}{{\mathsf{H}}}
\newcommand{\ii}{{\underline{i}}}
\newcommand{\jj}{{\underline{j}}}
\newcommand{\la}{{\lambda}}
\newcommand{\m}{{\mathfrak{m}}}

\newcommand{\ot}{{\otimes}}
\newcommand{\otot}{{\otimes\cdots\otimes}}
\newcommand{\too}{{\longrightarrow}}

\begin{document}

\title{Reflection functors and symplectic reflection algebras
for wreath products}

\author{Wee Liang Gan}
\address{Department of Mathematics,  Massachusetts Institute
of Technology, Cambridge, MA 02139, USA}
\email{wlgan@math.mit.edu}

\begin{abstract}
We construct reflection functors on categories of
modules over deformed wreath products of the
preprojective algebra of a quiver.
These functors give equivalences of categories 
associated to generic parameters which are 
in the same orbit under the Weyl group action.
We give applications to the representation theory of
symplectic reflection algebras of wreath product groups.
\end{abstract}


\maketitle

\section{\bf Introduction}

Deformed preprojective algebras $\Pi_\la$ associated to a quiver
were introduced by Crawley-Boevey and Holland  
in \cite{CBH}. A useful tool in their work is the reflection
functors, which gives an equivalence from the category of
modules over $\Pi_\la$ to the category of modules over
$\Pi_\mu$, when the parameters $\la$ and $\mu$ are in the
same orbit under an action of the Weyl group of the quiver
on the space of parameters.
Recently, in \cite{GG}, a one-parameter deformation
$\A_{n, \la, \nu}$ of the wreath product of
$\Pi_\la$ with $S_n$ was constructed.
The purpose of this paper is to generalize the construction of the
reflection functors to the algebras $\A_{n, \la, \nu}$.

Actually, we construct reflection functors $F_i$ for
the simple reflections $s_i$ at vertices $i$ without edge-loop. 
The author does not know if compositions of the functors 
satisfy the Weyl group relations.
However, we prove in Theorem \ref{equivalence} 
that if $\la, \nu$ are generic, then $F_i^2 \cong 1$. 

It is interesting that
to each $\A_{n, \la, \nu}$-module $V$, there is
a natural complex $\CC^\bullet(V)$, depending on $i$,
with the property that $F_i(V)=H^0(\CC^\bullet(V))$.
Assuming that $\la$ is generic and $\nu=0$, we prove
$H^r(\CC^\bullet(V))=0$ for all $r>0$, and hence
obtain a `dimension vector' formula for $F_i(V)$.

When the quiver is affine Dynkin of type ADE, there is a
finite subgroup $\Gamma \subset SL_2(\C)$ associated to it
by the McKay correspondence. 
Let $\GG$ be the wreath product group $S_n\ltimes\Gamma^n$.
In \cite{EG}, Etingof and Ginzburg introduced the 
symplectic reflection algebras $\hh_{t,k,c}(\GG)$
attached to $\GG$. A Morita equivalence between the algebras 
$\A_{n, \la, \nu}$ and $\hh_{t,k,c}(\GG)$ was constructed 
in \cite{GG}; in the case $n=1$, this was done in \cite{CBH}.
As a consequence, we obtain reflection functors for the algebras
$\hh_{t,k,c}(\GG)$. 

This paper is a step towards the classification of the symplectic 
reflection algebras of wreath product groups up to
Morita equivalence.
The reflection functors are defined when $\Gamma\neq\{1\}$.
Let us mention that when $\Gamma=\{1\}$ or $\Z/2\Z$,
the algebras $\hh_{t,k,c}(\GG)$ are the rational Cherednik 
algebras of type A or B. 
In these cases, Morita equivalences for the algebras
were constructed in \cite{BEG1} using shift functors,
and a complete classification in type A (which corresponds
to the affine Dynkin quiver of type $\mathrm A_0$)
was proved in \cite{BEG2} (for generic parameter).

We are interested in the representation theory of 
$\hh_{t,k,c}(\GG)$ when the parameter $t$ is nonzero.
When $n=1$, there is no parameter $k$, and
the finite dimensional simple modules of $\hh_{t,c}(\Gamma)$
were classified in \cite{CBH}.
When $n>1$, we have $\hh_{t,0,c}(\GG)
=\hh_{t,c}(\Gamma)^{\ot n}\rtimes \C[S_n]$.
Thus, there is also a classification of finite dimensional simple 
modules of $\hh_{t,0,c}(\GG)$.

In \cite{M2} (which generalizes \cite{EM}), Montarani found 
sufficient conditions for the existence of a deformation of
a finite dimensional simple $\hh_{1,0,c_0}(\GG)$-module to a 
$\hh_{1,k,c_0+c}(\GG)$-module for formal parameters $k,c$.
The proofs in \cite{EM} and \cite{M2} are
based on homological arguments.
We shall give a new proof by constructing the deformation
using reflection functors.
Moreover, using the reflection functors,
we will show that if a finite dimensional simple
$\hh_{1,0,c_0}(\GG)$-module can be formally deformed to
a $\hh_{1,k,c_0+c}(\GG)$-module, then
the conditions in \cite{M2} must necessarily hold.
We shall also use the reflection functors to prove 
the existence of certain flat families of finite 
dimensional simple $\hh_{t,k,c}(\GG)$-modules
(for complex parameters).

We expect that there will be other applications of
the reflection functors.

This paper is organized as follows.
In Section 2, we will recall the definition of the
algebras $\A_{n,\la,\nu}$, and construct 
the reflection functors. 
In Section 3 and Section 4, we give the proofs of
several identities required in the construction of the
reflection functors. 
In Section 5, we prove that the reflection functor
is an equivalence of categories for generic parameters.
We also construct the complex $\CC^\bullet(V)$, and
prove some other properties of the reflection functors. 
In Section 6, we give the applications
to the symplectic reflection algebras for wreath products.

\section{\bf Construction of the reflection functors} 

\subsection{}
We first recall some standard notions.
Let $\k$ be a commutative ring with $1$.
We shall work over $\k$. 

Let $Q$ be a quiver, and denote by
$I$ the set of vertices of $Q$.
The double $\QQ$ of $Q$ is the quiver obtained from $Q$ by
adding a reverse edge $\stackrel{a^*}{j\to i}$ 
for each edge $\stackrel{a}{i\to j}$ in $Q$.
We let $(a^*)^*:=a$ for any edge $a\in Q$. 
If $\stackrel{a}{i\to j}$ is an edge
in $\QQ$, we call $t(a):=i$ its tail, and 
$h(a):=j$ its head.
When $t(a)=h(a)$, we say that $a$ is an edge-loop.

The Ringel form of $Q$ is the bilinear form on $\Z^I$
defined by
$$\langle \al,\beta \rangle
:= \sum_{i\in I} \al_i\beta_i
- \sum_{a\in Q} \al_{t(a)}\beta_{h(a)}, \quad\mbox{ where }
\alpha=(\al_i)_{i\in I},
\ \beta=(\beta_i)_{i\in I}.$$
Let $(\al,\beta):= \langle \al,\beta \rangle
+\langle \beta,\al \rangle$ be its symmetrization.
We write $\epsilon_i\in \Z^I$ for the coordinate vector
corresponding to the vertex $i\in I$.
If there is no edge-loop at the vertex $i$, then
there is a reflection $s_i:\Z^I\to\Z^I$ defined by
$s_i(\al):= \al - (\al,\epsilon_i)\epsilon_i$.
We call $s_i$ a simple reflection.
The Weyl group $W$ is the group of automorphisms of $\Z^I$
generated by all the simple reflections.

Let $B:=\bigoplus_{i\in I} \k$,
and $E$ the free $\k$-module with
basis formed by the set of edges $\{a\in \QQ\}$.
Thus, $E$ is naturally a $B$-bimodule and $E=\bigoplus_{i,j\in I} 
E_{i,j}$,
where $E_{i,j}$ is spanned by the edges $a\in\QQ$ with 
$h(a)=i$
and $t(a)=j$. The path algebra of $\QQ$
is $\k\QQ := T_B E = \bigoplus_{n\geq 0} T^n_B E$, where
$T^n_B E = E\ot_B \cdots \ot_B E$ is the $n$-fold
tensor product.
The trivial path for the vertex $i$ is denoted by
$e_i$, an idempotent in $B$.
For any element $\la\in B$,
we will write $\la= \sum_{i\in I}\la_i e_i$ where 
$\la_i\in\k$.
If $w\in W$, $\la\in B$ and $\al\in \Z^I$, then
$(w\la)\cdot\al := \la\cdot(w^{-1}\al)$.
The reflection $r_i:B\to B$ dual to $s_i$ is defined
by $(r_i\la)_j := \la_j-(\epsilon_i,\epsilon_j)\la_i$.

\subsection{}
In this subsection, we recall the definition of
the algebra $\A_{n,\la,\nu}$ from \cite[Definition 1.2.3]{GG}.

From now on, we fix a positive integer $n$. 
Denote by $S_n$ the permutation group of 
$[1,n] := \{1, \ldots, n\}$,
and write $s_{ij}\in S_n$ for the transposition 
$i\leftrightarrow j$. Let $\B := B^{\ot n}$. 
For any $\ell\in [1,n]$, define the $\B$-bimodules
$$ \E_{\ell}:= B^{\ot (\ell-1)}\ot E\ot B^{\ot (n-\ell)}
\qquad \mathrm{and}\qquad \E :=
\bigoplus_{1\leq\ell\leq n} \E_\ell\,.$$
Given $\ell\in [1,n]$, $a\in \k\QQ$, 
and $\ii=(i_1, \ldots, i_n)\in I^n$, 
we write 
$$a_\ell\big|_{\ii} 
\quad\mbox{ for the element }\quad
e_{i_1}\otot ae_{i_\ell}\otot e_{i_n} \in T_\B\E_{\ell}.$$
We shall simply write $\big|_{\ii}$
for the element $e_{i_1}\otot e_{i_n}$.
If $a\in \QQ$ and $i_\ell=t(a)$, then let 
$$a_\ell(\ii):= (i'_1,\ldots,i'_n)\in I^n,\quad \mbox{ where }\quad
i'_m = \left\{ \begin{array}{ll}
i_m & \mbox{ if } m\neq\ell,\\
h(a) & \mbox{ if } m=\ell. \end{array} \right.$$

\begin{definition} \label{algebra}
For any $\la\in B$ and $\nu\in \k$, define the $\B$-algebra
$\A_{n,\la,\nu}$ to be the quotient 
of $T_{\B}\E\rtimes \k[S_n]$
by the following relations.
\begin{itemize}
\item[{\rm (i)}]
For any $\ell\in [1,n]$ and $\ii=(i_1,\ldots,i_n)\in I^n$:
$$
\Big(\sum_{a\in Q}
[a, a^*]-\la \Big)_\ell\Big|_{\ii}
= \nu \sum_{\{ m\neq\ell \mid i_m=i_\ell\}}
s_{\ell m}\Big|_{\ii}.
$$
\item[{\rm (ii)}]
For any $\ell,m\in [1,n]$ ($\ell\neq m$),
$a,b\in\QQ$, and $\ii=(i_1,\ldots,i_n)\in I^n$
with $i_\ell=t(a)$, $i_m=t(b)$:
$$
a_\ell \big|_{b_m(\ii)} 
b_m\big|_{\ii}
-b_m\big|_{a_\ell(\ii)} 
a_\ell\big|_{\ii}  
 =  \left\{ \begin{array}{ll}
\nu s_{\ell m} \big|_{\ii}
& \textrm{if $b\in Q$ and $a=b^*$},\\
- \nu s_{\ell m} \big|_{\ii}
& \textrm{if $a\in Q$ and $b=a^*$},\\
0 & \textrm{else}\,.  \end{array} \right.  
$$
\end{itemize}
\end{definition}

If $n=1$, there is no parameter $\nu$,
and $\A_{1,\la}$ is the deformed preprojective
algebra $\Pi_\la$. Observe that 
$\A_{n,\la,0} = \Pi_\la^{\ot n}\rtimes \k[S_n]$.

We shall denote by $\A_{n,\la,\nu}\mathrm{-mod}$
the category of left $\A_{n,\la,\nu}$-modules.

\subsection{}
Let $i$ be a vertex of $Q$ such that 
there is no edge-loop at $i$.
We shall define the reflection functor 
$$ F_i : \A_{n,\la,\nu}\mathrm{-mod}
\ \too\ \A_{n,r_i\la,\nu}\mathrm{-mod}. $$
In the case $n=1$, the reflection functors were constructed by
Crawley-Boevey and Holland in \cite[Theorem 5.1]{CBH}.
They were also constructed by Nakajima
in the context of quiver varieties, see \cite[Remark 3.20]{Na}.
They are similar to (but not the same as) 
the reflection functors of 
Bernstein, Gelfand and Ponomarev in \cite{BGP}.

Let $V$ be a $\A_{n,\la,\nu}$-module.
We will first define $F_i(V)$ as a $\B\rtimes \k[S_n]$-module.

Up to isomorphism,
the algebra $\A_{n,\la,\nu}$ does not depend on
the orientation of $Q$, so we may 
assume that $i$ is a sink in $Q$; 
let $$R:=\{a\in Q\mid h(a)=i\}.$$
For any $\jj=(j_1,\ldots,j_n)\in I^n$, let 
$$V_{\jj} := \big|_\jj V, \quad\mbox{ and }\quad
\Delta(\jj):=\{m\in [1,n]\mid j_m=i\}.$$
For any $D\subseteq \Delta(\jj)$, define
$$\X (D) := \mbox{ the set of all maps }\xi:D\to R:
m\mapsto \xi(m).$$ 
Given $\xi\in \X(D)$, let
$$t(\jj,\xi) := (t_1,\ldots,t_n)\in I^n, 
\quad \mbox{ where }\quad
t_m = \left\{ \begin{array}{ll}
j_m & \mbox{ if } m\notin D,\\
t(\xi(m)) & \mbox{ if } m\in D. \end{array} \right.$$
Define 
$$V(\jj,D) := 
\bigoplus_{\xi\in \X(D)} V_{t(\jj,\xi)}.$$
In particular, $V(\jj,\emptyset) = V_\jj$. Write 
$$\pi_{\jj,\xi}:V(\jj,D)\too V_{t(\jj,\xi)},\qquad
\mu_{\jj,\xi}:V_{t(\jj,\xi)}\too V(\jj,D)$$ 
for the projection map and inclusion map, respectively.

If $\sigma\in S_n$, then let
$\sigma(\jj) := (j_{\sigma^{-1}(1)}, \ldots, j_{\sigma^{-1}(n)})$.
We have $\Delta(\sigma(\jj)) = \sigma(\Delta(\jj))$. 
If $\xi\in \X(D)$, then let $\sigma(\xi) \in \X(\sigma(D))$
be the map $m\mapsto \xi(\sigma^{-1}m)$. Let
$$ \sigma\big|_\jj: V(\jj, D) \too V(\sigma(\jj), \sigma(D)),
\qquad \sigma\big|_\jj := \sum_{\xi\in \X(D)}
\mu_{\sigma(\jj), \sigma(\xi)} \sigma \pi_{\jj,\xi}.$$

Suppose $p\in D$. We have a restriction map
$\rho_p: \X(D)\to \X(D\setminus\{p\})$.
For each $\xi\in \X(D)$, we have a composition of maps
$$ \xymatrix{V(\jj,D) \ar[rr]^{\pi_{\jj,\xi}} & & 
V_{t(\jj,\xi)} \ar[rr]^{\xi(p)_p\big|_{t(\jj,\xi)}} & & 
V_{t(\jj,\rho_p(\xi))} \ar@{^{(}->}[rr]^{\mu_{\jj,\rho_p(\xi)}}
& & V(\jj, D\setminus\{p\})}. $$
We also have a composition of maps
$$ \xymatrix{ V(\jj, D\setminus\{p\}) \ar[rr]^{\pi_{\jj,\rho_p(\xi)}}
& & V_{t(\jj,\rho_p(\xi))} 
       \ar[rr]^{\xi(p)^*_p\big|_{t(\jj,\rho_p(\xi))}} & & 
V_{t(\jj,\xi)} \ar@{^{(}->}[rr]^{\mu_{\jj,\xi}} & & V(\jj,D) }.$$
Define 
$$\pi_{\jj,p}: V(\jj,D)\too V(\jj,D\setminus\{p\}), \quad
 \pi_{\jj,p} := \sum_{\xi\in \X(D)}
\mu_{\jj,\rho_p(\xi)} \xi(p)_p\big|_{t(\jj,\xi)} \pi_{\jj,\xi},$$
and
$$\mu_{\jj,p}: V(\jj,D\setminus\{p\})\too V(\jj,D), \quad
\mu_{\jj,p}:= \sum_{\xi\in \X(D)}
\mu_{\jj,\xi} \xi(p)^*_p\big|_{t(\jj,\rho_p(\xi))}
\pi_{\jj,\rho_p(\xi)}.$$
The maps $\pi_{\jj,p}$ and $\mu_{\jj,p}$ depend on
$D$, but we suppress it from our notations.

We state here the following lemma which will be 
used later.
\begin{lemma} \label{firstlemma}
For any $p,q\in D$ ($p\neq q$), we have the following.

{\rm (i)}  $\pi_{\sigma(\jj), \sigma(p)} \sigma\big|_\jj
= \sigma\big|_{\jj} \pi_{\jj,p}$, and
$\mu_{\sigma(\jj), \sigma(p)} \sigma\big|_\jj
= \sigma\big|_{\jj} \mu_{\jj,p}$.

{\rm (ii)} $\pi_{\jj,p}\mu_{\jj,p} = \lambda_i + 
\nu \sum_{m\in \Delta(\jj)\setminus D}s_{pm}\big|_\jj$.

{\rm (iii)} $\pi_{\jj,p}\mu_{\jj,q}
= \mu_{\jj,q}\pi_{\jj,p} -\nu s_{pq}\big|_{\jj}$.

{\rm (iv)} $\pi_{\jj,p}\pi_{\jj,q}=\pi_{\jj,q}\pi_{\jj,p}$,
and $\mu_{\jj,p}\mu_{\jj,q}=\mu_{\jj,q}\mu_{\jj,p}$.
\end{lemma}
The proof of Lemma \ref{firstlemma} will be given in Section 3.
Let
$$ V_\jj(D):=\left\{ \begin{array}{ll}
\bigcap_{p\in D} \mathrm{Ker}(\pi_{\jj,p}) & \mbox{ if }
D \neq\emptyset,\\
V_\jj & \mbox{ if } D=\emptyset. \end{array} \right.$$
Let $V'_\jj := V_\jj(\Delta(\jj))$.

\begin{definition}
Let $F_i(V) := V'= \bigoplus_{\jj\in I^n} V'_\jj$ as a 
$\B\rtimes \k[S_n]$-module.
\end{definition}

\subsection{}
The proofs of Lemmas \ref{lemmaforcase1}, \ref{lemmaforcase2}
and \ref{lemmaforcase3} of this subsection will be given
in Section 3.
The proof of Proposition \ref{functor} of this subsection
will be given in Section 4.

Given any $\ell\in[1,n]$, $a\in \QQ$, 
and $\jj=(j_1,\ldots, j_n)\in I^n$ with $j_\ell=t(a)$,
we shall define a map
$a'_\ell\big|_\jj : V'_\jj \to V'_{a_\ell(\jj)}$.
There are three cases.

{\bf Case (I)}, $h(a),t(a)\neq i$: 
Then $\ell\notin\Delta(\jj)=\Delta(a_\ell(\jj))$. 

Let $D\subset \Delta(\jj)$.
For any $\xi\in \X(D)$, we have a composition of maps
$$ \xymatrix{ V(\jj,D) \ar[rr]^{\pi_{\jj,\xi}} &&
V_{t(\jj,\xi)} \ar[rr]^{a_\ell\big|_{t(\jj,\xi)}} 
&& V_{t(a_\ell(\jj),\xi)}
\ar@{^{(}->}[rr]^{\mu_{a_\ell(\jj),\xi}} && 
V\big(a_\ell(\jj),D)} .$$ 
Define 
$$a_\ell\big|_{\jj,D}: V(\jj,D)\too 
V\big(a_\ell(\jj),D),\qquad 
a_\ell\big|_{\jj,D}:= \sum_{\xi\in \X(D)}
\mu_{a_\ell(\jj),\xi} a_\ell\big|_{t(\jj,\xi)} \pi_{\jj,\xi}.$$
Let $a'_\ell\big|_\jj := a_\ell\big|_{\jj,\Delta(\jj)}$.

\begin{lemma} \label{lemmaforcase1}
{\rm (i)} If $p\in D$, then $\pi_{a_\ell(\jj),p} a_\ell\big|_{\jj,D}
= a_\ell\big|_{\jj,D\setminus\{p\}} \pi_{\jj,p}$.

{\rm (ii)} If $p\notin D$, then
$\mu_{a_\ell(\jj),p} a_\ell\big|_{\jj,D}
= a_\ell\big|_{\jj,D\cup\{p\}} \mu_{\jj,p}$.
\end{lemma}

It follows from Lemma \ref{lemmaforcase1}(i) that
$a'_\ell\big|_\jj$ defined in Case (I) sends 
$V'_\jj$ into $V'_{a_\ell(\jj)}$.

{\bf Case (II)}, $t(a)=i$: Then $\ell\in\Delta(\jj)$,
and $\Delta(a_\ell(\jj)) = \Delta(\jj)\setminus\{\ell\}$.

Suppose $\ell\in D\subset \Delta(\jj)$.
For any $r\in R$, we have an injective map
$$ \tau_{r,\ell,D}: \X(D\setminus\{\ell\}) \hookrightarrow \X(D):\eta
\mapsto \tau_{r,\ell,D}(\eta), $$
where 
$$ \tau_{r,\ell,D}(\eta)(m) := \left\{ \begin{array}{ll}
\eta(m) & \mbox{ if } m\in D\setminus\{\ell\},\\
r & \mbox{ if } m=\ell. \end{array} \right.$$
Since $t(\jj,\tau_{r,\ell,D}(\eta)) = t(r^*_\ell(\jj), \eta)$, there is
a projection map
$$ \tau_{r,\ell,\jj,D}^!: V(\jj,D) \too 
V(r^*_\ell(\jj),D\setminus\{\ell\}),
\quad \tau_{r,\ell,\jj,D}^!:= \sum_{\eta\in\X(D\setminus\{\ell\})}
\mu_{r^*_\ell(\jj), \eta} \pi_{\jj,\tau_{r,\ell,D}(\eta)},$$
and an inclusion map 
$${\tau_{r,\ell,\jj,D}}_!:
V(r^*_\ell(\jj),D\setminus\{\ell\}) \too V(\jj,D), \quad
{\tau_{r,\ell,\jj,D}}_!:= \sum_{\eta\in\X(D\setminus\{\ell\})}
\mu_{\jj,\tau_{r,\ell,D}(\eta)}\pi_{r^*_\ell(\jj), \eta}.$$
Let $a'_\ell\big|_\jj := 
\tau_{a^*, \ell, \jj, \Delta(\jj)}^!$.

\begin{lemma} \label{lemmaforcase2}
{\rm (i)} $\sum_{r\in R} 
{\tau_{r,\ell,\jj,D}}_! \tau_{r,\ell,\jj,D}^! =1$.

{\rm (ii)} 
If $p\in D\setminus\{\ell\}$, then 
$$\pi_{r^*_\ell(\jj),p} \tau_{r,\ell,\jj,D}^!
= \tau_{r, \ell, \jj, D\setminus\{p\}}^! \pi_{\jj,p}
\quad \mbox{ and } \quad
\pi_{\jj,p} {\tau_{r,\ell,\jj,D}}_!
= {\tau_{ r,\ell,\jj,D\setminus\{p\} }}_! \pi_{r^*_\ell(\jj),p}.$$

{\rm (iii)}
If $p\notin D$, then
$$ \mu_{r^*_\ell(\jj),p} \tau^!_{r,\ell,\jj,D}
= \tau^!_{r,\ell,\jj,D\cup\{p\}} \mu_{\jj,p}
\quad \mbox{ and } \quad
\mu_{\jj,p} { \tau_{r,\ell,\jj,D} }_! 
= { \tau_{r,\ell,\jj,D\cup\{p\}} }_! 
 \mu_{r^*_\ell(\jj),p}.$$
\end{lemma}

It follows from Lemma \ref{lemmaforcase2}(ii) that
$a'_\ell\big|_\jj$ defined in Case (II) sends
$V'_\jj$ into $V'_{a_\ell(\jj)}$.

{\bf Case (III)}, $h(a)=i$:
Then $\ell\notin\Delta(\jj)$, and $\Delta(a_\ell(\jj))=
\Delta(\jj)\cup\{\ell\}$.

Let $D\subset \Delta(\jj)$. We have the inclusion map
$$ { \tau_{a, \ell, a_\ell(\jj), D\cup\{\ell\}} }_! :
V(\jj,D) \too V(a_\ell(\jj), D\cup\{\ell\}).$$
Define
$$ \theta_{a,\ell,\jj, D} : V(\jj,D) \too
V(a_\ell(\jj), D\cup\{\ell\})$$
by
$$ \theta_{a,\ell,\jj, D} :=
\Big(-\lambda_i + \mu_{a_\ell(\jj), \ell} \pi_{a_\ell(\jj), \ell}
  + \nu\sum_{m\in D} s_{m\ell}\big|_{a_\ell(\jj)} \Big) 
  { \tau_{a, \ell, a_\ell(\jj), D\cup\{\ell\}} }_!.$$
Let $a'_\ell\big|_\jj :=  \theta_{a,\ell,\jj, \Delta(\jj)}$. 

\begin{lemma} \label{lemmaforcase3}
{\rm (i)} If $p\in D$, then
$\pi_{a_\ell(\jj),p} \theta_{a,\ell,\jj, D}
= \theta_{a,\ell,\jj, D\setminus\{p\}} \pi_{\jj,p}$.

{\rm (ii)} We have 
$$ \pi_{a_\ell(\jj), \ell} \theta_{a,\ell,\jj, \Delta(\jj)}
= \nu\sum_{m\in \Delta(\jj)} s_{m\ell} \big|_{a_\ell(\jj)} 
   {\tau_{ a, \ell, a_\ell(\jj), \Delta(a_\ell(\jj))\setminus\{m\} }}_!
   \pi_{\jj,m} .$$

{\rm (iii)} If $p\in D$, then
$ \theta_{a,\ell,\jj, D} \mu_{\jj,p}
= \mu_{a_\ell(\jj), p} \theta_{a,\ell,\jj, D\setminus\{p\}}$.
\end{lemma}

It follows from Lemma \ref{lemmaforcase3}(i)-(ii) that
$a'_\ell\big|_\jj$ defined in Case (III) maps
$V'_\jj$ into $V'_{a_\ell(\jj)}$.

Thus, $F_i(V)$ is a $T_\B \E \rtimes \k[S_n]$-module,
where $a_\ell\big|_\jj\in \E$ acts by 
$a'_\ell\big|_\jj$.

\begin{proposition} \label{functor}
With the above action,
$F_i(V)$ is a $\A_{n, r_i \la, \nu}$-module.
\end{proposition}

It is clear from our construction
that the assignment $V\mapsto F_i(V)$ is functorial.

\section{\bf Proofs of lemmas}

{\bf Proof of Lemma \ref{firstlemma}}:

(i) \begin{align*}
\pi_{\sigma(\jj), \sigma(p)} \sigma\big|_\jj
=& \sum_{\xi\in \X(D)}
\mu_{\sigma(\jj), \rho_{\sigma(p)} (\sigma(\xi))}
\sigma(\xi)(\sigma(p))_{\sigma(p)}
  \big|_{t(\sigma(\jj), \sigma(\xi))}
\sigma \pi_{\jj,\xi}\\
=& \sum_{\xi\in \X(D)}
\mu_{\sigma(\jj), \rho_{\sigma(p)} (\sigma(\xi))}
\sigma
\xi(p)_p\big|_{t(\jj,\xi)} \pi_{\jj,\xi}\\
=& \sigma\big|_{\jj} \mu_{\jj,p}.
\end{align*}

\begin{align*}
\mu_{\sigma(\jj), \sigma(p)} \sigma\big|_\jj
=& \sum_{\xi\in\X(D)}
\mu_{\sigma(\jj),\sigma(\xi)}
\sigma(\xi)(\sigma(p))^*_{\sigma(p)}
  \big|_{t(\sigma(\jj), \rho_{\sigma(p)}(\sigma(\xi))}
  \sigma \pi_{\jj, \rho_p(\xi)}\\
=& \sum_{\xi\in\X(D)}
\mu_{\sigma(\jj),\sigma(\xi)}
\sigma
\xi(p)^*_p\big|_{t(\jj,\rho_p(\xi))}
\pi_{\jj, \rho_p(\xi)}\\
=& \sigma\big|_{\jj} \mu_{\jj,p}.
\end{align*}

(ii) \begin{align*}
\pi_{\jj,p}\mu_{\jj,p} 
=& \sum_{\xi\in \X(D)} \mu_{\jj,\rho_p(\xi)}
  \xi(p)_p \big|_{t(\jj,\xi)} \xi(p)^*_p \big|_{t(\jj,\rho_p(\xi))}
  \pi_{\jj,\rho_p(\xi)} \\
=& \sum_{\eta\in \X(D\setminus\{p\})}
  \mu_{\jj,\eta}\Big( \sum_{a\in R} aa^* \Big)_p\Big|_{t(\jj,\eta)}
  \pi_{\jj,\eta} \\
=& \sum_{\eta\in \X(D\setminus\{p\})}
  \mu_{\jj,\eta}\Big( \lambda_i 
   + \nu \sum_{m\in \Delta(\jj)\setminus D} s_{pm} \Big)
  \pi_{\jj,\eta} \\
=& \lambda_i+\nu\sum_{m\in \Delta(\jj)\setminus D}s_{pm}\big|_\jj.
\end{align*}

(iii) \begin{align*}
\pi_{\jj,p}\mu_{\jj,q} 
=& \sum_{\xi\in \X(D)} \mu_{\jj,\rho_p(\xi)}
  \xi(p)_p \big|_{t(\jj,\xi)} \xi(q)^*_q \big|_{t(\jj,\rho_q(\xi)}
  \pi_{\jj,\rho_q(\xi)} \\
=& \sum_{\xi\in \X(D)} \mu_{\jj,\rho_p(\xi)}
   \xi(q)^*_q \big|_{t( \jj,\rho_p(\rho_q(\xi)))}
  \xi(p)_p \big|_{t(\jj,\rho_q(\xi))} 
  \pi_{\jj,\rho_q(\xi)} \\
 & - \nu \sum_{ \{ \xi\in\X(D) \mid \xi(p)=\xi(q)  \} }
  \mu_{\jj,\rho_p(\xi)} s_{pq}
  \pi_{\jj,\rho_q(\xi)} \\
=& \mu_{\jj,q}\pi_{\jj,p} -\nu s_{pq}\big|_\jj.  
\end{align*}

(iv) \begin{align*}
\pi_{\jj,p} \pi_{\jj,q}
=& \sum_{\xi \in \X(D)}  
  \mu_{\jj, \rho_p(\rho_q(\xi))}
  \xi(p)_p \big|_{ t(\jj, \rho_q(\xi)) }
  \xi(q)_q \big|_{ t(\jj,\xi) }
  \pi_{\jj,\xi} \\
=& \sum_{\xi \in \X(D)}
    \mu_{\jj, \rho_p(\rho_q(\xi))}
    \xi(q)_q \big|_{ t(\jj, \rho_p(\xi)) }
    \xi(p)_p \big|_{ t(\jj,\xi) }
    \pi_{\jj,\xi} \\
=& \pi_{\jj,q} \pi_{\jj,p}.
\end{align*}

\begin{align*}
\mu_{\jj,p} \mu_{\jj,q}
=& \sum_{\xi\in\X(D)}
  \mu_{\jj,\xi} \xi(p)^*_p \big|_{t(\jj,\rho_p(\xi))}
  \xi(q)^*_q\big|_{t(\jj, \rho_q(\rho_p(\xi)))}
  \pi_{\jj, \rho_q(\rho_p(\xi))}\\
=& \sum_{\xi\in\X(D)}
  \mu_{\jj,\xi} \xi(q)^*_q \big|_{t(\jj,\rho_q(\xi))}
  \xi(p)^*_p\big|_{t(\jj, \rho_q(\rho_p(\xi)))}
  \pi_{\jj, \rho_q(\rho_p(\xi))}\\
=& \mu_{\jj,q} \mu_{\jj,p}.
\end{align*} \qed

{\bf Proof of Lemma \ref{lemmaforcase1}}:

(i) \begin{align*}
\pi_{a_\ell(\jj),p} a_\ell\big|_{\jj,D} 
=& \sum_{\xi\in\X(D)} 
   \mu_{a_\ell(\jj),\rho_p(\xi)}
   \xi(p)_p\big|_{t(a_\ell(\jj),\xi)}
   a_\ell\big|_{t(\jj,\xi)} 
   \pi_{\jj,\xi}\\
=& \sum_{\xi\in\X(D)}
   \mu_{a_\ell(\jj),\rho_p(\xi)}
   a_\ell\big|_{t(a_\ell(\jj),\rho_p(\xi))}
   \xi(p)_p\big|_{t(\jj,\xi)}
   \pi_{\jj,\xi}\\
=& a_\ell\big|_{\jj,D\setminus\{p\}} \pi_{\jj,p}.
\end{align*}

(ii) \begin{align*}
\mu_{a_\ell(\jj),p} a_\ell\big|_{\jj,D}
=& \sum_{\xi\in \X(D\cup\{p\})}
\mu_{a_\ell(\jj),\xi} \xi(p)^*_p
 \big|_{t(a_\ell(\jj), \rho_p(\xi))}
 a_\ell\big|_{t(\jj,\rho_p(\xi))} \pi_{\jj,\rho_p(\xi)} \\
=& \sum_{\xi\in \X(D\cup\{p\})}
\mu_{a_\ell(\jj),\xi} a_\ell\big|_{t(\jj,\xi)}
 \xi(p)^*_p \big|_{t(\jj,\rho_p(\xi))}
 \pi_{\jj,\rho_p(\xi)} \\
=& a_\ell\big|_{\jj,D\cup\{p\}} \mu_{\jj,p}. 
\end{align*} \qed

{\bf Proof of Lemma \ref{lemmaforcase2}}:

(i) \begin{align*}
\sum_{r\in R} 
{\tau_{r,\ell,\jj,D}}_! \tau_{r,\ell,\jj,D}^!
=& \sum_{r\in R}  \sum_{\eta\in \X(D\setminus\{\ell\})}
\mu_{\jj, \tau_{r,\ell,D} (\eta)}
\pi_{\jj, \tau_{r,\ell,D} (\eta)} \\
=& \sum_{ \xi \in \X(D) }
\mu_{\jj,\xi} \pi_{\jj,\xi} \\
=& 1.
\end{align*}

(ii) \begin{align*}
\pi_{r^*_\ell(\jj),p} \tau_{r,\ell,\jj,D}^! 
=& \sum_{\eta\in \X(D\setminus\{\ell\})}
 \mu_{r^*_\ell(\jj),\rho_p(\eta)}
 \eta(p)_p\big|_{t(r^*_\ell(\jj),\eta)}
 \pi_{\jj,\tau_{r,\ell,D}(\eta)} \\
=& \sum_{ \{ \xi\in\X(D)\mid \xi(\ell)=r \} }
  \mu_{r^*_\ell(\jj),\rho_p(\rho_\ell(\xi))}
  \xi(p)_p\big|_{t(\jj,\xi)} 
  \pi_{\jj,\xi} \\
=& \Big( \sum_{ \{ \varepsilon \in \X(D\setminus\{p\}) \mid
\varepsilon(\ell)=r \} }
  \mu_{r^*_\ell(\jj),\rho_\ell(\varepsilon)}
  \pi_{\jj,\varepsilon} \Big) \pi_{\jj,p} \\
=& \Big( \sum_{\zeta\in \X(D\setminus\{ \ell,p \})}
  \mu_{r^*_\ell(\jj),\zeta} 
  \pi_{\jj, \tau_{r,\ell,D\setminus\{p\}} (\zeta)} \Big) 
  \pi_{\jj,p} \\
=& \tau_{r, \ell,\jj, D\setminus\{p\}}^! \pi_{\jj,p}.
\end{align*}

\begin{align*}
\pi_{\jj,p} {\tau_{r,\ell,\jj,D}}_!
=& \Big(\sum_{\xi\in\X(D)} 
 \mu_{\jj,\rho_p(\xi)} \xi(p)_p\big|_{t(\jj,\xi)} 
 \pi_{\jj,\xi} \Big)
 \Big( \sum_{\eta\in \X(D\setminus\{\ell\})}
 \mu_{\jj, \tau_{r,\ell,D} (\eta)} \pi_{r^*_\ell(\jj),\eta} \Big) \\
=& \sum_{\eta\in \X(D\setminus\{\ell\})}
   \mu_{\jj,\rho_p(\tau_{r,\ell,D} (\eta))}
  \eta(p)_p \big|_{t( \jj, \tau_{r,\ell,D} (\eta) )} 
  \pi_{r^*_\ell(\jj),\eta} \\
=& \Big( \sum_{\zeta\in \X(D\setminus \{\ell,p\} )}
   \mu_{\jj, \tau_{r, \ell,D\setminus\{p\}} (\zeta)  } 
  \pi_{r^*_\ell(\jj),\zeta} \Big)
  \pi_{r^*_\ell(\jj),p} \\
=& { \tau_{r, \ell, \jj, D\setminus \{ p\}} }_! \pi_{r^*_\ell(\jj),p}
\end{align*}

(iii) \begin{align*}
\mu_{r^*_\ell(\jj),p} \tau^!_{r,\ell,\jj,D}
=& \sum_{\xi\in \X(D\cup\{p\}\setminus\{\ell\})}
\mu_{r^*_\ell(\jj), \xi} \xi(p)^*_p
 \big|_{t(r^*_\ell(\jj), \rho_p(\xi))}
\pi_{\jj, \tau_{r,\ell,D}(\rho_p(\xi))} \\
=& \sum_{ \{\zeta\in\X(D\cup\{p\}) \mid \zeta(\ell)=r  \} }
\mu_{r^*_\ell(\jj),\rho_\ell(\zeta)} 
\zeta(p)^*_p \big|_{t(r^*_\ell(\jj),\rho_p(\rho_\ell(\zeta)))}
\pi_{\jj, \rho_p(\zeta)} \\
=& \Big( \sum_{\eta\in \X(D\cup\{p\}\setminus\{\ell\})}
\mu_{r^*_\ell(\jj),\eta} \pi_{\jj,\tau_{r,\ell,D\cup\{p\}}(\eta)}
  \Big)  \mu_{\jj,p} \\
=& \tau^!_{r,\ell,\jj,D\cup\{p\}} \mu_{\jj,p}.
\end{align*}

\begin{align*}
\mu_{\jj,p} { \tau_{r,\ell,\jj,D} }_!
=& \sum_{\{ \xi\in\X(D\cup\{p\}) \mid \xi(\ell)=r \}}
\mu_{\jj,\xi} \xi(p)^*_p \big|_{t(\jj,\rho_p(\xi))}
\pi_{r^*_\ell(\jj), \rho_p(\rho_\ell(\xi))} \\
=& \sum_{\eta\in \X(D\cup\{p\}) \setminus \{\ell\}}
\mu_{\jj,\tau_{r,\ell,D\cup\{p\}}(\eta)}
\eta(p)^*_p \big|_{ t(r^*_\ell(\jj), \rho_p(\eta)) }
\pi_{r^*_\ell(\jj), \rho_p(\eta)} \\
=& { \tau_{r,\ell,\jj,D\cup\{p\}} }_!
 \mu_{r^*_\ell(\jj),p}.
\end{align*} \qed

{\bf Proof of Lemma \ref{lemmaforcase3}}:

(i) Using Lemma \ref{firstlemma}(iii), Lemma \ref{firstlemma}(iv), 
and Lemma \ref{lemmaforcase2}(ii),
we have
\begin{align*}
\pi_{a_\ell(\jj),p} \theta_{a,\ell,\jj, D}
=& \Big( -\lambda_i \pi_{a_\ell(\jj),p}
 + \mu_{a_{\ell(\jj),\ell}} \pi_{a_\ell(\jj), \ell}
    \pi_{a_\ell(\jj), p}
 -\nu s_{p\ell} \big|_{a_\ell(\jj)}  
  \pi_{a_\ell(\jj), \ell} \\
& + \nu \sum_{m\in D} \pi_{a_\ell(\jj), p} 
   s_{m\ell}\big|_{a_\ell(\jj)} \Big)
  { \tau_{a,  \ell, a_\ell(\jj), D\cup\{\ell\}} }_! \\
=& -\lambda_i 
   { \tau_{a,  \ell, a_\ell(\jj), D\cup\{\ell\}\setminus\{p\}} }_!
     \pi_{\jj,p} 
   + \mu_{a_{\ell(\jj),\ell}} \pi_{a_\ell(\jj), \ell}
    { \tau_{a,  \ell, a_\ell(\jj), D\cup\{\ell\}\setminus\{p\}} }_!
    \pi_{\jj,p} \\
&  + \nu \sum_{m\in D\setminus\{p\}}
   s_{m\ell} \big|_{a_\ell(\jj)} \pi_{a_\ell(\jj), p}
   { \tau_{a,  \ell, a_\ell(\jj), D\cup\{\ell\}} }_! \\
=& \Big( -\lambda_i  
   + \mu_{a_{\ell(\jj),\ell}} \pi_{a_\ell(\jj), \ell}
   + \nu \sum_{m\in D\setminus\{p\}} 
     s_{m\ell} \big|_{a_\ell(\jj)}  \Big)
  { \tau_{a,  \ell, a_\ell(\jj), D\cup\{\ell\}\setminus\{p\}} }_!
  \pi_{\jj,p} \\
=& \theta_{a,\ell,\jj, D\setminus\{p\}} \pi_{\jj,p}.
\end{align*}

(ii) Using Lemma \ref{firstlemma}(ii)
and Lemma \ref{lemmaforcase2}(ii), we have
\begin{align*}
\pi_{a_\ell(\jj),\ell} \theta_{a,\ell,\jj, \Delta(\jj)}
=& \Big( -\lambda_i \pi_{a_\ell(\jj),\ell} 
+ \lambda_i \pi_{a_\ell(\jj),\ell}
+ \nu \sum_{m\in \Delta(\jj)} 
  s_{m\ell} \big|_{a_\ell(\jj)} \pi_{a_\ell(\jj),m} \Big)    
  { \tau_{a,  \ell, a_\ell(\jj), \Delta(\jj)\cup\{\ell\}} }_! \\
=& \nu \sum_{m\in \Delta(\jj)}
  s_{m\ell} \big|_{a_\ell(\jj)}
    { \tau_{a,  \ell, a_\ell(\jj), 
              \Delta(\jj)\cup\{\ell\}\setminus\{m\}} }_!
     \pi_{\jj,m} . 
\end{align*} 

(iii) Using Lemma \ref{lemmaforcase2}(iii)
and Lemma \ref{firstlemma}(iii), we have
\begin{align*}
\theta_{a,\ell,\jj,D}\mu_{\jj,p}
=& \Big(-\la_i + \mu_{a_\ell(\jj),\ell}\pi_{a_\ell(\jj),\ell}
+\nu \sum_{m\in D} s_{m\ell} \big|_{a_\ell(\jj)}  \Big)
{ \tau_{a,\ell,a_\ell(\jj),D\cup\{\ell\}} }_!
\mu_{\jj,p} \\
=&  \Big( -\la_i + \mu_{a_\ell(\jj),\ell}\pi_{a_\ell(\jj),\ell}
+\nu \sum_{m\in D} s_{m\ell} \big|_{a_\ell(\jj)}  \Big)
\mu_{a_\ell(\jj),p}
{ \tau_{a,\ell,a_\ell(\jj),D\cup\{\ell\}\setminus\{p\}} }_! \\
=&   \Big( -\la_i \mu_{a_\ell(\jj),p}
+ \mu_{a_\ell(\jj),\ell}\big(
        \mu_{a_\ell(\jj),p} \pi_{a_\ell(\jj),\ell}
        - \nu s_{p\ell} \big|_{a_\ell(\jj)} \big) \\
 & +\nu \sum_{m\in D} s_{m\ell} \big|_{a_\ell(\jj)}  
\mu_{a_\ell(\jj),p}\Big)
{ \tau_{a,\ell,a_\ell(\jj),D\cup\{\ell\}\setminus\{p\}} }_! \\
=& \mu_{a_\ell(\jj),p} \theta_{a,\ell,\jj,D\setminus\{p\}} .
\end{align*} \qed

\section{\bf Proof of Proposition \ref{functor}}

To prove Proposition \ref{functor},
we have to check that the relations in Definition 
\ref{algebra} hold for the maps $a'_\ell\big|_\jj$.

First, we check the relations of type (i).
Suppose we are given $\ell\in [1,n]$
and $\jj=(j_1,\ldots,j_n)\in I^n$.
There are two cases.

{\bf Case 1}, $j_\ell=i$: Then $\ell\in\Delta(\jj)$.
We have
\begin{align*}
\sum_{a\in R} a'_\ell \big|_{a^*_\ell(\jj)}
  a^{*'}_\ell \big|_\jj
=& \sum_{a\in R} 
  \theta_{a,\ell,a^*_\ell(\jj),\Delta(\jj)\setminus\{\ell\}}
  \tau_{a,\ell,\jj,\Delta(\jj)}^! \\
=& \sum_{a\in R}
  \Big( -\la_i + \mu_{\jj,\ell} \pi_{\jj,\ell}
       + \nu\sum_{m\in \Delta(\jj)\setminus\{\ell\}}
                 s_{m\ell} \big|_\jj \Big) 
  { \tau_{a,\ell,\jj,\Delta(\jj)} }_!
  \tau_{a,\ell,\jj,\Delta(\jj)}^! \\
=& -\la_i + \mu_{\jj,\ell} \pi_{\jj,\ell}
       + \nu\sum_{m\in \Delta(\jj)\setminus\{\ell\}}
                 s_{m\ell} \big|_\jj 
   \quad \mbox{ (by Lemma \ref{lemmaforcase2}(i)) } \\
=& (r_i\la)_i + \nu\sum_{m\in \Delta(\jj)\setminus\{\ell\}}
                 s_{m\ell} \big|_\jj
   \quad \mbox{ (since $\pi_{\jj,\ell}$ vanishes on $V'_\jj$) }.
\end{align*}

{\bf Case 2}, $j_\ell\neq i$: Then $\ell\notin \Delta(\jj)$.
If $a\in R$ and $t(a)=j_\ell$, then
\begin{align*}
a^{*'}_\ell \big|_{ a_\ell(\jj) } a'_\ell \big|_\jj
=& \tau^!_{a,\ell, a_\ell(\jj), \Delta(\jj)\cup\{\ell\}}
   \theta_{ a,\ell,\jj,\Delta(\jj) } \\
=& \Big( \sum_{\xi\in\X(\Delta(\jj))}
     \mu_{\jj,\xi} \pi_{a_\ell(\jj), 
             \tau_{a,\ell, \Delta(\jj)\cup\{\ell\}} (\xi)} \Big)\\
 & \times
  \Big( -\la_i + \mu_{a_\ell(\jj),\ell} \pi_{a_\ell(\jj),\ell}
     + \nu\sum_{m\in \Delta(\jj)} s_{m\ell}\big|_{a_\ell(\jj)}
   \Big) \\
 & \times
  \Big(\sum_{\varepsilon\in\X(\Delta(\jj))}
     \mu_{a_\ell(\jj), 
              \tau_{a,\ell, \Delta(\jj)\cup\{\ell\}} (\varepsilon)}  
     \pi_{\jj,\varepsilon} \Big) \\
=& -\la_i 
   + \sum_{\xi\in \X(\Delta(\jj))} 
    \mu_{\jj,\xi} (a^*a)_\ell\big|_{t(\jj,\xi)} \pi_{\jj,\xi} \\
  & + \nu \sum_{m\in\Delta(\jj)}
    \sum_{ \{ \xi\in\X(\Delta(\jj)) \mid \xi(m)=a \} } 
   \mu_{\jj,\xi} s_{m\ell} \pi_{\jj,\xi} .
\end{align*}
Hence,
\begin{align*}
& \sum_{ \{a\in Q\mid h(a)=j_\ell\} }  
  a'_\ell\big|_{a^*_\ell(\jj)} a^{*'}_\ell\big|_\jj
- \sum_{ \{a\in Q\mid t(a)=j_\ell\} }
  a^{*'}_\ell\big|_{a_\ell(\jj)} a'_\ell\big|_\jj \\
=& \sum_{\xi\in \X(\Delta(\jj))}
 \mu_{\jj,\xi}   
  \Big( \sum_{ \{a\in Q\mid h(a)=j_\ell\} } aa^*
       -\sum_{ \{a\in Q\mid t(a)=j_\ell\} } a^*a
  \Big)_\ell \big|_{t(\jj,\xi)}
 \pi_{\jj,\xi} \\
& + \la_i \Big(\sum_{ \{a\in R\mid t(a)=j_\ell\} } 1\Big)
   - \nu \sum_{ \{a\in R\mid t(a)=j_\ell\} }
    \sum_{m\in\Delta(\jj)} 
    \sum_{ \{ \xi\in\X(\Delta(\jj)) \mid \xi(m)=a \} }
   \mu_{\jj,\xi} s_{m\ell} \pi_{\jj,\xi} \\
=& \sum_{\xi\in \X(\Delta(\jj))}
   \mu_{\jj,\xi}
  \Big( \la_{j_\ell} 
  + \nu\sum_{ \{m\neq\ell \mid j_m=j_\ell\} }
          s_{m\ell} \big|_{t(\jj,\xi)}
  + \nu\sum_{ \{ m\in\Delta(\jj) \mid t(\xi(m))=j_\ell \} } 
          s_{m\ell} \big|_{t(\jj,\xi)}
  \Big) \pi_{\jj,\xi} \\
& -(\epsilon_{j_\ell},\epsilon_i)\la_i
     - \nu \sum_{ \{a\in R\mid t(a)=j_\ell\} }
    \sum_{m\in\Delta(\jj)}
    \sum_{ \{ \xi\in\X(\Delta(\jj)) \mid \xi(m)=a \} }
   \mu_{\jj,\xi} s_{m\ell} \pi_{\jj,\xi} \\
=& (r_i\la)_{j_\ell}
  + \nu\sum_{ \{m\neq\ell \mid j_m=j_\ell\} }
          s_{m\ell} \big|_\jj .
\end{align*}

Next, we check the relations of type (ii).
Suppose we are given $\ell, m\in [1,n]$ ($\ell\neq m$),
$a,b\in \QQ$, and $\jj=(j_1,\ldots,j_n)\in I^n$,
such that $j_\ell=t(a)$, $j_m=t(b)$.
There are six cases.

{\bf Case 1}, $t(a),h(a),t(b),h(b)\neq i$:
\begin{align*}
& a'_\ell\big|_{b_m(\jj)} b'_m\big|_\jj 
- b'_m\big|_{a_\ell(\jj)} a'_\ell\big|_\jj \\
=& a_\ell\big|_{b_m(\jj), \Delta(\jj)} b_m\big|_{\jj,\Delta(\jj)} 
- b_m\big|_{a_\ell(\jj),\Delta(\jj)} a_\ell\big|_{\jj,\Delta(\jj)}\\
=& \sum_{\xi\in \X(\Delta(\jj))}
\mu_{a_\ell(b_m(\jj)), \xi}
\Big( a_\ell\big|_{t(b_m(\jj),\xi)} b_m\big|_{t(\jj,\xi)}
  - b_m\big|_{t(a_\ell(\jj),\xi)} a_\ell\big|_{t(\jj,\xi)} \Big)
\pi_{\jj,\xi} \\
=& \left\{ \begin{array}{ll}
\nu s_{\ell m} \big|_{\jj}
& \textrm{if $b\in Q$ and $a=b^*$},\\
- \nu s_{\ell m} \big|_{\jj}
& \textrm{if $a\in Q$ and $b=a^*$},\\
0 & \textrm{else}\,.  \end{array} \right.  
\end{align*}

{\bf Case 2}, $t(a)=i$ and $t(b),h(b)\neq i$:
\begin{align*}
& a'_\ell\big|_{b_m(\jj)} b'_m\big|_\jj 
- b'_m\big|_{a_\ell(\jj)} a'_\ell\big|_\jj \\
=& \tau^!_{a^*, \ell, b_m(\jj), \Delta(\jj)}
b_m\big|_{\jj, \Delta(\jj)}
- b_m\big|_{a_\ell(\jj), \Delta(\jj)\setminus\{\ell\}}
\tau^!_{a^*,\ell,\jj,\Delta(\jj)}\\
=& \sum_{\eta\in\X(\Delta(\jj)\setminus\{\ell\})}
  \mu_{a_\ell(b_m(\jj)),\eta}
  b_m\big|_{t(\jj, \tau_{a^*,\ell,\Delta(\jj)}(\eta) )}
  \pi_{\jj, \tau_{a^*,\ell,\Delta(\jj)}(\eta)} \\
 & - \sum_{\eta\in\X(\Delta(\jj)\setminus\{\ell\})}
  \mu_{b_m(a_\ell(\jj)),\eta}
  b_m\big|_{t(a_\ell(\jj),\eta)}
  \pi_{\jj, \tau_{a^*,\ell,\Delta(\jj)}(\eta)} \\
=& 0.
\end{align*}

{\bf Case 3}, $h(a)=i$ and $t(b),h(b)\neq i$:
We have 
\begin{align*}
{ \tau_{a, \ell, a_\ell(b_m(\jj)), \Delta(\jj)\cup\{\ell\}} }_!
b_m\big|_{\jj,\Delta(\jj)}
=& \sum_{\xi\in\X(\Delta(\jj)}
\mu_{ a_\ell(b_m(\jj)),\tau_{a,\ell,\Delta(\jj)\cup\{\ell\}} }
b_m\big|_{t(\jj,\xi)} \pi_{\jj,\xi}\\
=& b_m\big|_{a_\ell(\jj), \Delta(\jj)\cup\{\ell\}}
  { \tau_{a,\ell,a_\ell(\jj),\Delta(\jj)\cup\{\ell\}} }_!.
\end{align*}
By Lemma \ref{lemmaforcase1},
\begin{align*}
\mu_{a_\ell(b_m(\jj)),\ell}
\pi_{a_\ell(b_m(\jj)),\ell}
b_m\big|_{a_\ell(\jj),\Delta(\jj)\cup\{\ell\}}
=& \mu_{a_\ell(b_m(\jj)),\ell}
   b_m\big|_{a_\ell(\jj),\Delta(\jj)}
   \pi_{a_\ell(\jj),\ell} \\
=& b_m\big|_{a_\ell(\jj),\Delta(\jj)\cup\{\ell\}}
   \mu_{a_\ell(\jj),\ell} 
   \pi_{a_\ell(\jj),\ell} .
\end{align*}
Hence,
\begin{align*}
& a'_\ell\big|_{b_m(\jj)} b'_m\big|_\jj 
- b'_m\big|_{a_\ell(\jj)} a'_\ell\big|_\jj \\
=& \theta_{a,\ell,b_m(\jj),\Delta(\jj)}
b_m\big|_{\jj,\Delta(\jj)}
- b_m\big|_{a_\ell(\jj),\Delta(\jj)\cup\{\ell\}}
\theta_{a,\ell,\jj,\Delta(\jj)} \\
=& \Big(  -\la_i
+ \mu_{a_\ell(b_m(\jj)),\ell} \pi_{a_\ell(b_m(\jj)),\ell}
+ \nu \sum_{p\in \Delta(\jj)}
           s_{p\ell} \big|_{a_\ell(b_m(\jj))} \Big)
 { \tau_{a, \ell, a_\ell(b_m(\jj)), \Delta(\jj)\cup\{\ell\}} }_!
b_m\big|_{\jj,\Delta(\jj)} \\
 & - b_m\big|_{a_\ell(\jj),\Delta(\jj)\cup\{\ell\}}
 \Big( -\la_i + \mu_{a_\ell(\jj),\ell} \pi_{a_\ell(\jj),\ell}
  + \nu \sum_{p\in \Delta(\jj)}
           s_{p\ell} \big|_{a_\ell(\jj)} \Big)
 { \tau_{a,\ell,a_\ell(\jj),\Delta(\jj)\cup\{\ell\}} }_! \\
=& 0.
\end{align*}

{\bf Case 4}, $t(a)=t(b)=i$:
\begin{align*}
a'_\ell\big|_{b_m(\jj)} b'_m\big|_\jj 
=& \tau^!_{a^*,\ell,b_m(\jj),\Delta(\jj)\setminus\{m\}}
\tau^!_{b^*,m,\jj,\Delta(\jj)} \\
=& \sum_{\eta\in\X( \Delta(\jj)\setminus \{\ell,m\} )}
  \mu_{a_\ell(b_m(\jj)), \eta}
  \pi_{ \jj, \tau_{b^*, m, \Delta(\jj)} 
              (\tau_{a^*,\ell,\Delta(\jj)\setminus\{m\}} (\eta)) }\\
=& \tau^!_{b^*,m,a_\ell(\jj),\Delta(\jj)\setminus\{\ell\}}
\tau^!_{a^*,\ell,\jj,\Delta(\jj)} \\
=& b'_m\big|_{a_\ell(\jj)} a'_\ell\big|_\jj .
\end{align*}

{\bf Case 5}, $h(a)=t(b)=i$: 
\begin{align*}
& b'_m\big|_{a_\ell(\jj)} a'_\ell\big|_\jj \\
=& \tau^!_{b^*,m,a_\ell(\jj), \Delta(\jj)\cup\{\ell\} }
   \theta_{a,\ell, \jj, \Delta(\jj)} \\
=& \tau^!_{b^*,m,a_\ell(\jj), \Delta(\jj)\cup\{\ell\} }
   \Big( -\la_i 
     + \mu_{a_\ell(\jj),\ell} \pi_{a_\ell(\jj),\ell}
         + \nu \sum_{p\in \Delta(\jj)}
           s_{p\ell} \big|_{a_\ell(\jj)} \Big)
 { \tau_{a,\ell,a_\ell(\jj),\Delta(\jj)\cup\{\ell\}} }_! \\
=&\Big( -\la_i
+ \mu_{a_\ell(b_m(\jj)),\ell} \pi_{a_\ell(b_m(\jj)),\ell}
+ \nu \sum_{p\in \Delta(\jj)\setminus\{m\}}
           s_{p\ell} \big|_{a_\ell(b_m(\jj))} \Big) 
    \tau^!_{b^*,m,a_\ell(\jj), \Delta(\jj)\cup\{\ell\} } \\ 
  &\times
    { \tau_{a,\ell,a_\ell(\jj),\Delta(\jj)\cup\{\ell\}} }_! 
  + \nu s_{m\ell} \big|_{b_\ell(a_\ell(\jj))}
    \tau^!_{b^*,\ell,a_\ell(\jj), \Delta(\jj)\cup\{\ell\} }
    { \tau_{a,\ell,a_\ell(\jj),\Delta(\jj)\cup\{\ell\}} }_! \\
& \mbox{ (using Lemma \ref{lemmaforcase2}).}      
\end{align*}
Now
\begin{align*}
& \tau^!_{b^*,m,a_\ell(\jj), \Delta(\jj)\cup\{\ell\} } 
    { \tau_{a,\ell,a_\ell(\jj),\Delta(\jj)\cup\{\ell\}} }_! \\
=& \Big( \sum_{\xi\in \X(\Delta(\jj)\cup\{\ell\}\setminus\{m\})}
  \mu_{b_m(a_\ell(\jj)),\xi}
  \pi_{a_\ell(\jj), 
      \tau_{ b^*,m,\Delta(\jj)\cup\{\ell\} } (\xi)} \Big)
   \Big( \sum_{\zeta\in \X(\Delta(\jj))}
      \mu_{ a_\ell(\jj), 
          \tau_{ a,\ell,\Delta(\jj)\cup\{\ell\} } (\zeta) }
      \pi_{\jj,\zeta} \Big)\\
=& \sum_{\eta \in \X( \Delta(\jj)\setminus\{m\} )}
   \mu_{ a_\ell(b_m(\jj)),
         \tau_{ a,\ell,\Delta(\jj)\cup\{\ell\}\setminus\{m\} } 
               (\eta) }
   \pi_{\jj, \tau_{b^*,m,\Delta(\jj)} (\eta)} \\
=& { \tau_{a,\ell,a_\ell(b_m(\jj)),
        \Delta(\jj)\cup\{\ell\}\setminus\{m\}} }_!
   \tau^!_{b^*,m,\jj,\Delta(\jj)},
\end{align*}
and
\begin{align*}
& \tau^!_{b^*,\ell,a_\ell(\jj), \Delta(\jj)\cup\{\ell\} }
  { \tau_{a,\ell,a_\ell(\jj),\Delta(\jj)\cup\{\ell\}} }_! \\
=&  \Big( \sum_{\xi\in \X(\Delta(\jj))}
  \mu_{b_\ell(a_\ell(\jj)),\xi}
  \pi_{a_\ell(\jj), 
      \tau_{ b^*,\ell,\Delta(\jj)\cup\{\ell\} } (\xi)} \Big)
   \Big( \sum_{\zeta\in \X(\Delta(\jj))}
      \mu_{ a_\ell(\jj), 
          \tau_{ a,\ell,\Delta(\jj)\cup\{\ell\} } (\zeta) }
      \pi_{\jj,\zeta} \Big)\\
=& \left\{ \begin{array}{ll}
   0 & \mbox{ if } a\neq b^*,\\
   1 & \mbox{ if } a = b^*.  \end{array}\right.
\end{align*}
Hence,
$$
b'_m\big|_{a_\ell(\jj)} a'_\ell\big|_\jj \\
= \left\{ \begin{array}{ll}
a'_\ell\big|_{b_m(\jj)} b'_m\big|_\jj 
 & \mbox{ if } a\neq b^*,\\
a'_\ell\big|_{b_m(\jj)} b'_m\big|_\jj 
 +\nu s_{m\ell}\big|_{\jj}
 & \mbox{ if } a=b^*.
\end{array}\right.
$$

{\bf Case 6}, $h(a)=h(b)=i$:
\begin{align*}
& a'_\ell\big|_{b_m(\jj)} b'_m\big|_\jj  \\
=& \theta_{a,\ell,b_m(\jj),\Delta(\jj)\cup\{m\}}
   \theta_{b,m,\jj,\Delta(\jj)}  \\
=& \Big( -\la_i 
   + \mu_{a_\ell(b_m(\jj)),\ell} \pi_{a_\ell(b_m(\jj)),\ell}
   + \nu\sum_{p \in \Delta(\jj)\cup\{m\}}
          s_{p\ell}\big|_{a_\ell(b_m(\jj))} \Big)
  { \tau_{a,\ell,a_\ell(b_m(\jj)),\Delta(\jj)\cup\{\ell,m\}} }_! \\
 & \times
  \Big( -\la_i
   + \mu_{b_m(\jj),m} \pi_{b_m(\jj),m}
   + \nu\sum_{q\in \Delta(\jj)} s_{qm}\big|_{b_m(\jj)} \Big)
  { \tau_{b,m,b_m(\jj),\Delta(\jj)\cup\{m\}} }_! \\
=& \Big( -\la_i 
   + \mu_{a_\ell(b_m(\jj)),\ell} \pi_{a_\ell(b_m(\jj)),\ell}
   + \nu\sum_{p \in \Delta(\jj)\cup\{m\}}
          s_{p\ell}\big|_{a_\ell(b_m(\jj))} \Big) \\
  & \times
  \Big( -\la_i 
   + \mu_{a_\ell(b_m(\jj)),m} \pi_{a_\ell(b_m(\jj)),m}
   + \nu\sum_{q \in \Delta(\jj)}
          s_{qm}\big|_{a_\ell(b_m(\jj))} \Big) \\
  & \times 
  { \tau_{a,\ell,a_\ell(b_m(\jj)),\Delta(\jj)\cup\{\ell,m\}} }_! 
  { \tau_{b,m,b_m(\jj),\Delta(\jj)\cup\{m\}} }_!
  \quad \mbox{ (using Lemma \ref{lemmaforcase2}) }.
\end{align*}
Now
\begin{align*}
&  { \tau_{a,\ell,a_\ell(b_m(\jj)),\Delta(\jj)\cup\{\ell,m\}} }_! 
  { \tau_{b,m,b_m(\jj),\Delta(\jj)\cup\{m\}} }_! \\
=& \sum_{\xi\in\X(\Delta(\jj))}
\mu_{a_\ell(b_m(\jj)), 
  \tau_{a,\ell, \Delta(\jj)\cup\{\ell,m\}}
      (\tau_{b,m,\Delta(\jj)\cup\{m\}} (\xi) )}
\pi_{\jj,\xi}\\
=& { \tau_{b,m,a_\ell(b_m(\jj)),\Delta(\jj)\cup\{\ell,m\}} }_!
  { \tau_{a,\ell,a_\ell(\jj),\Delta(\jj)\cup\{\ell\}} }_! 
\end{align*}
Using Lemma \ref{firstlemma}, we have
\begin{align*}
& \big(\mu_{a_\ell(b_m(\jj)),\ell} \pi_{a_\ell(b_m(\jj)),\ell}  
+\nu s_{m\ell}\big|_{a_\ell(b_m(\jj))} \big)
\mu_{a_\ell(b_m(\jj)),m} \pi_{a_\ell(b_m(\jj)),m} \\
=& \mu_{a_\ell(b_m(\jj)),\ell} \mu_{a_\ell(b_m(\jj)),m} 
\pi_{a_\ell(b_m(\jj)),\ell}  \pi_{a_\ell(b_m(\jj)),m} \\
=& \big(\mu_{a_\ell(b_m(\jj)),m} \pi_{a_\ell(b_m(\jj)),m}  
+\nu s_{m\ell}\big|_{a_\ell(b_m(\jj))} \big)
\mu_{a_\ell(b_m(\jj)),\ell} \pi_{a_\ell(b_m(\jj)),\ell}.
\end{align*}
Moreover,
\begin{align*}
\sum_{p\in\Delta(\jj)\cup\{m\}} s_{p\ell}
\sum_{q \in \Delta(\jj)} s_{qm}
=& \sum_{\substack{p,q\in\Delta(\jj) \\ p\neq q}} s_{p\ell}s_{qm} 
+ \sum_{q\in\Delta(\jj)} s_{q\ell}s_{qm}
+ \sum_{q\in\Delta(\jj)} s_{m\ell}s_{qm}\\
=& \sum_{\substack{p,q\in\Delta(\jj) \\ p\neq q}} s_{qm}s_{p\ell} 
+ \sum_{q\in\Delta(\jj)} s_{m\ell}s_{q\ell}
+ \sum_{q\in\Delta(\jj)} s_{qm}s_{q\ell}\\
=& \sum_{p\in\Delta(\jj)\cup\{\ell\}} s_{pm}
\sum_{q \in \Delta(\jj)} s_{q\ell}.
\end{align*}
Hence,
$$ a'_\ell\big|_{b_m(\jj)} b'_m\big|_\jj 
= b'_m\big|_{a_\ell(\jj)} a'_\ell\big|_\jj .$$

This completes the proof of Proposition \ref{functor}.

\section{\bf Properties of the reflection functors}

\subsection{}
Let $i$ be a vertex of $Q$ such that there is no 
edge-loop at $i$.
Define $\Lambda_i$ to be the set of all $(\la,\nu)\in B\times \k$
such that $\la_i \pm \nu \sum_{m=2}^r s_{1 m}$
are invertible in $\k[S_r]$ for all $r\in [1,n]$.

\begin{theorem} \label{equivalence}
If $(\la,\nu)\in\Lambda_i$, then
the functor $$F_i: \A_{n,\la,\nu}\mathrm{-mod}
\ \too\ \A_{n,r_i\la,\nu}\mathrm{-mod}$$ is an equivalence
of categories with quasi-inverse functor $F_i$.
\end{theorem}

\begin{proof}
Let $V\in\A_{n,\la,\nu}\mathrm{-mod}$ and
$V'=F_i(V)\in\A_{n,r_i\la,\nu}\mathrm{-mod}$.
Let $\jj\in I^n$.
                                                                          
By Lemma \ref{firstlemma}(ii),
for any $p\in \Delta(\jj)$, the composition
$$ \xymatrix{
V(\jj, \Delta(\jj)\setminus\{p\})
\ar[rr]^{\mu_{\jj,p}} &&
V(\jj,\Delta(\jj))
\ar[rr]^{\pi_{\jj,p}} &&
V(\jj, \Delta(\jj)\setminus\{p\})
}$$ is equal to $\la_i$.
Since $\la_i$ is invertible, we have a direct sum decomposition
$$ V(\jj,\Delta(\jj)) = \mathrm{Ker} \pi_{\jj,p}
\oplus \mathrm{Im} \mu_{\jj,p}.$$
                                                                          
Now $V'(\jj,\Delta(\jj))=V(\jj,\Delta(\jj))$, and
$$
V'(\jj,\Delta(\jj)\setminus\{p\})
= \bigoplus_{\eta\in\X(\Delta(\jj)\setminus\{p\})}
  V'_{t(\jj,\eta)}
\subset
\bigoplus_{\eta\in\X(\Delta(\jj)\setminus\{p\})}
  \Big( \bigoplus_{r\in R}
   V_{t(\jj, \tau_{r,p,\Delta(\jj)}(\eta))}
  \Big)
= V(\jj,\Delta(\jj)).
$$
Observe that the kernel of the map
$$ -\la_i + \mu_{\jj,p}\pi_{\jj,p}:V(\jj,\Delta(\jj)) \too
V(\jj,\Delta(\jj)) $$
is $\mathrm{Im} \mu_{\jj,p} \subset V(\jj,\Delta(\jj))$.
Hence,
$$ F_i(F_i(V))_\jj = F_i(V')_\jj = \left\{ \begin{array}{ll}
\bigcap_{p\in \Delta(\jj)} \mathrm{Im}(\mu_{\jj,p}) & \mbox{ if }
\Delta(\jj) \neq\emptyset,\\
V_\jj & \mbox{ if } \Delta(\jj)=\emptyset. \end{array} \right.$$
                                                                          
Suppose $\Delta(\jj) = \{ p_1, \ldots, p_r \}$.
We have a canonical map
$$\mu_{\jj,p_1} \cdots \mu_{\jj,p_r}:V_\jj \too V(\jj,\Delta(\jj)).$$
By Lemma \ref{firstlemma}(iv), this map does not depend on
the ordering of $p_1,\ldots,p_r$, and its image lies in
$\bigcap_{p\in \Delta(\jj)} \mathrm{Im}(\mu_{\jj,p})$.
We claim that it is an isomorphism from $V_\jj$ to
$\bigcap_{p\in \Delta(\jj)} \mathrm{Im}(\mu_{\jj,p})$.
By Lemma  \ref{firstlemma}(ii), each
$\mu_{\jj,p}$ is injective. Hence, we have to show that
$\bigcap_{p\in \Delta(\jj)} \mathrm{Im}(\mu_{\jj,p})
\subset \mathrm{Im}(\mu_{\jj,p_1} \cdots \mu_{\jj,p_r})$.
                                                                          
It suffices to prove that, for $h\in [2,r]$, we have
$$ \mathrm{Im}(\mu_{\jj,p_1} \cdots \mu_{\jj,p_{h-1}})
\cap \mathrm{Im}(\mu_{\jj,p_h}) \subset
\mathrm{Im}(\mu_{\jj,p_1} \cdots \mu_{\jj,p_h})
\subset V(\jj,\Delta(\jj)).$$
Suppose that
$\mu_{\jj,p_1} \cdots \mu_{\jj,p_{h-1}}(v)
= \mu_{\jj,p_h}(w)$.
Then, by Lemma \ref{firstlemma},
\begin{gather*}
\la_i w = \pi_{\jj,p_h}\mu_{\jj,p_h}(w)
=  \pi_{\jj,p_h} \mu_{\jj,p_1} \cdots \mu_{\jj,p_{h-1}}(v) \\
=  \mu_{\jj,p_1} \cdots \mu_{\jj,p_{h-1}} \pi_{\jj,p_h} (v)
  -\nu\sum_{g=1}^{h-1} s_{p_g p_h}
    \mu_{\jj,p_1} \cdots \mu_{\jj,p_{g-1}}\mu_{\jj,p_{g+1}}
    \cdots \mu_{\jj,p_{h-1}} (v).
\end{gather*}
Hence,
\begin{align*}
\la_i \mu_{\jj,p_h}(w)
=& \mu_{\jj,p_1} \cdots \mu_{\jj,p_{h}} \pi_{\jj,p_h} (v)
 -\nu\sum_{g=1}^{h-1} s_{p_g p_h}
  \mu_{\jj,p_1} \cdots  \mu_{\jj,p_{g}}  \cdots \mu_{\jj,p_{h-1}} 
(v)\\
=& \mu_{\jj,p_1} \cdots \mu_{\jj,p_{h}} \pi_{\jj,p_h} (v)
 -\nu\sum_{g=1}^{h-1} s_{p_g p_h} \mu_{\jj,p_h}(w).
\end{align*}
It follows that
$$ \mu_{\jj,p_h}(w)
= (\la_i+\nu\sum_{g=1}^{h-1} s_{p_g p_h})^{-1}
  \mu_{\jj,p_1} \cdots \mu_{\jj,p_{h}} \pi_{\jj,p_h} (v)
\in \mathrm{Im}(\mu_{\jj,p_1} \cdots \mu_{\jj,p_h}).$$
                                                                          
It remains to show that our map $V\to F_i(F_i(V))$
commutes with the actions. Let
$a_\ell\big|_\jj\in\E$, with $t(a)=j_\ell$.
Let $\Delta(\jj)=\{p_1,\ldots,p_r\}$.
                                                                          
If $h(a),t(a)\neq i$, then by Lemma \ref{lemmaforcase1}(ii),
$$ a_\ell\big|_{\jj, \Delta(\jj)}
\mu_{\jj,p_1} \cdots \mu_{\jj,p_r}
= \mu_{a_\ell(\jj),p_1} \cdots \mu_{a_\ell(\jj),p_r}
a_\ell\big|_\jj.$$
                                                                          
If $t(a)=i$, let $p_r=\ell$.
Then by Lemma \ref{lemmaforcase2}(iii),
\begin{align*}
\tau^!_{a^*,\ell,\jj,\Delta(\jj)}
\mu_{\jj,p_1} \cdots \mu_{\jj,p_r}
=& \mu_{a_\ell(\jj),p_1} \cdots \mu_{a_\ell(\jj),p_{r-1}}
\tau^!_{a^*,\ell,\jj,\{\ell\}} \mu_{\jj,p_r}\\
=& \mu_{a_\ell(\jj),p_1} \cdots \mu_{a_\ell(\jj),p_{r-1}}
a_\ell\big|_\jj.
\end{align*}
                                                                          
If $h(a)=i$, then the map
$(F_iF_iV)_\jj\to (F_iF_iV)_{a_\ell(\jj)}$ is
\begin{gather*}
\Big( \la_i
 +\big( -\la_i+ \mu_{a_\ell(\jj),\ell} \pi_{a_\ell(\jj),\ell} \big)
 + \nu \sum_{m\in \Delta(\jj)} s_{m\ell}\big|_{a_\ell(\jj)} \Big)
 { \tau_{a,\ell,a_\ell(\jj),\Delta(\jj)\cup\{\ell\}} }_! \\
= \la_i { \tau_{a,\ell,a_\ell(\jj),D\cup\{\ell\}} }_!
+\theta_{a,\ell,\jj,D}.
\end{gather*}
By Lemma \ref{lemmaforcase2}(iii) and
Lemma \ref{lemmaforcase3}(iii), we have
\begin{align*}
& \big(\la_i { \tau_{a,\ell,a_\ell(\jj),\Delta(\jj)\cup\{\ell\}} }_!
+\theta_{a,\ell,\jj,\Delta(\jj)}\big)
\mu_{\jj,p_1} \cdots \mu_{\jj,p_r} \\
=&  \mu_{a_\ell(\jj),p_1} \cdots \mu_{a_\ell(\jj),p_r}
\big( \la_i { \tau_{a,\ell,a_\ell(\jj), \{\ell\}} }_!
 + \theta_{a,\ell,\jj,\emptyset} \big) \\
=& \mu_{a_\ell(\jj),p_1} \cdots \mu_{a_\ell(\jj),p_r}
 \mu_{a_\ell(\jj),\ell} \pi_{a_\ell(\jj),\ell}
 { \tau_{a,\ell,a_\ell(\jj),\{\ell\}} }_! \\
=& \mu_{a_\ell(\jj),p_1} \cdots \mu_{a_\ell(\jj),p_r}
  \mu_{a_\ell(\jj),\ell} a_\ell\big|_\jj. 
\end{align*}
\end{proof}
                                                                          
It is easy to see that the functor
$F_i: \A_{n,\la,\nu}\mathrm{-mod}
\,\to\,\A_{n,r_i\la,\nu}\mathrm{-mod}$
is left exact.

\begin{corollary} \label{exact}
Suppose $(\la,\nu)\in \Lambda_i$.
Then the functor $F_i:\A_{n,\la,\nu}\mathrm{-mod}
\to \A_{n,r_i\la,\nu}\mathrm{-mod}$ is exact.
\end{corollary}
\begin{proof}
By Theorem \ref{equivalence}, $F_i$ has a quasi-inverse
functor. Since $F_i$ is left exact, it must also be right
exact by \cite[Theorem 5.8.3]{We}. 
\end{proof}

Let $\k'$ be a commutative $\k$-algebra.
The following corollary 
will be used later in the proof of Proposition \ref{flat2}.

\begin{corollary} \label{fiber}
Suppose $(\la,\nu)\in \Lambda_i$.
Let $V$ be any $\A_{n,\la,\nu}$-module.
Then $F_i(V)\ot_\k \k'=F_i(V\ot_\k \k')$.
\end{corollary}
\begin{proof}
There is a natural map 
$f:F_i(V)\ot_\k \k' \to F_i(V\ot_\k \k')$.
We have
$$
(F_i(F_iV))\ot_\k \k' \too F_i((F_iV)\ot_\k \k') 
\too F_i(F_i(V\ot_\k \k')).
$$
By Theorem \ref{equivalence}, the composition of
these two maps is the identity map of $V\ot_\k \k'$.
The injectivity of the first map for any $V$
implies that $f$ is injective.
The surjectivity of the second map and the exactness
of $F_i$ imply that $f$ is surjective.
\end{proof}

\begin{remark}
Let $g$ be an automorphism of the graph $Q$.
Then $g$ acts on $B$ by $(g\la)_j=\la_{g^{-1}(j)}$ for any $j$,
and $g(r_i\la)=r_{g(i)}(g\la)$.
It was pointed out to the author by Iain Gordon that
$g$ induces an isomorphism of algebras
$\A_{n,\la,\nu}\to\A_{n,g\la,\nu}$.
Observe that the following diagram commutes:
$$\xymatrix{
\A_{n,g\la,\nu}-\mathrm{mod}
\ar[rr]^{}
\ar[d]_{F_{g(i)}} && \A_{n,\la,\nu}-\mathrm{mod}
\ar[d]^{F_i} \\
\A_{n,r_{g(i)}(g\la),\nu}-\mathrm{mod}
\ar[rr]^{}
&& \A_{n,r_i\la,\nu}-\mathrm{mod} }$$
\end{remark}

\subsection{}
In this subsection, we recall the definitions of a
commutative cube and its associated complex; 
see \cite[\S3]{Kh} for a more detailed discussion.

Let $\Delta$ be a finite set. For any $J\subset\Delta$,
we let $\Z J$ be the $\Z$-module freely 
generated by the elements of $J$, and write
$\det(J)$ for $\det(\Z J)$.
If $p\in \Delta\setminus J$, then we define an isomorphism
$$\iota: \det(J) \too \det(J\cup\{p\}),\quad
x\mapsto x\wedge p.$$

Let $\mathscr Z$ be the category of modules over a ring.

\begin{definition}
A commutative $\Delta$-cube $(Z,\psi)$ (over $\mathscr Z$)
consists of data:
\begin{itemize}
\item  an object $Z(J)\in \mathrm{Ob}(\mathscr Z)$
for each $J\subset\Delta$;
\item  a morphism $\psi_{J,p}: Z(J)\to Z(J\cup\{p\})$
for each $J\subset\Delta$ and $p\in \Delta\setminus J$.
\end{itemize}
These data are required to satisfy the following conditions:
for each $J\subset\Delta$ and $p,q\in \Delta\setminus J$
where $p\neq q$, we have
$\psi_{J\cup\{p\},q}\psi_{J,p}
=\psi_{J\cup\{q\},p}\psi_{J,q}$.
\end{definition}

Let $(Z,\psi)$ be a commutative $\Delta$-cube. 
We shall construct a complex $\CC^\bullet(Z)$
in the category $\mathscr Z$.
For each integer $r$, let
$$\CC^r(Z):= \bigoplus_{\substack{J\subset\Delta\\ |J|=r}} 
             Z(J) \ot_\Z \det(J).$$
Define the map
$$ d: \CC^r(Z)\to\CC^{r+1}(Z), \quad d:= 
\sum_{\substack{J\subset\Delta\\ |J|=r}}
\sum_{p\in \Delta\setminus J} \psi_{J,p} \ot \iota.$$
It is easy to check that $d^2=0$, so that
$(\CC^\bullet(Z), d)$ is a complex.

Let $q\in\Delta$, and let $\Delta_q:=\Delta\setminus\{q\}$.
Define
$$ Z_0(J):=Z(J),\quad Z_1(J):=Z(J\cup\{q\}),
\quad \mbox{ for all } J\subset\Delta_q.$$
Then both $(Z_0, \psi)$ and $(Z_1, \psi)$ are
commutative $\Delta_q$-cubes.
Let 
$$f: \CC^\bullet (Z_0) \too \CC^\bullet (Z_1),
\quad f := 
\sum_{J\subset\Delta_q} \psi_{J,q} \ot \mathrm{Id}.$$
The map $f$ is a morphism of complexes.
We note that the complex
$\CC^\bullet (Z)$ is the cone of the morphism $f[-1]$,
and we have a short exact sequence
\begin{equation} \label{cone}
0\too \CC^\bullet (Z_1)[-1] \too \CC^\bullet (Z)
\too \CC^\bullet (Z_0) \too 0.
\end{equation}

\begin{example} \label{idempotent}
Let $Z'$ be an object of $\mathscr Z$, and
let $\psi_1, \ldots, \psi_m$ be a set of
commuting endomorphisms of $Z'$.
Let $\Delta:=[1,m]$. For each $J\subset\Delta$, define
$$ Z(J):= \left\{\begin{array}{ll}
\mathrm{Im}(\psi_{q_1}\cdots\psi_{q_r}) &\mbox{ if }
  J = \{ q_1, \ldots, q_r \}, \\
Z' & \mbox{ if } J=\emptyset.
\end{array}\right. $$
If $J\subset\Delta$ and $p\in\Delta\setminus J$, then 
define $\psi_{J,p} : Z(J)\to Z(J\cup\{p\})$ to be the 
restriction of the morphism $\psi_p$ to $Z(J)$.
It is clear that $(Z,\psi)$ is a commutative $\Delta$-cube.

{\it Claim}: 
If $\psi_1, \ldots, \psi_m$ are idempotents,
then $H^r (\CC^\bullet(Z)) = 0$ for all $r>0$.

{\it Proof of Claim}:
This is clear if $m=1$. We shall prove the claim
by induction on $m$.

Let $q=1$. We have the commutative $\Delta_q$-cubes
$(Z_0,\psi)$ and $(Z_1,\psi)$. 
By (\ref{cone}), we have the long exact sequence
\begin{gather*}
0 \too H^0(\CC^\bullet (Z)) \too H^0(\CC^\bullet (Z_0))
\stackrel{\psi_1}{\too} H^0(\CC^\bullet (Z_1)) \too 
H^1(\CC^\bullet (Z)) \too \\
\too H^1(\CC^\bullet (Z_0)) \too H^1(\CC^\bullet (Z_1)) \too 
H^2(\CC^\bullet (Z)) \too H^2(\CC^\bullet (Z_0)) \too \cdots
\end{gather*}
By the induction hypothesis, we have
$H^r(\CC^\bullet (Z_0))=H^r(\CC^\bullet (Z_1))=0$
for all $r>0$. 
Hence, $H^r(\CC^\bullet (Z))=0$ for all $r>1$, and
$H^1(\CC^\bullet (Z))$ is isomorphic to the cokernel
of $\psi_1: H^0(\CC^\bullet (Z_0))\to H^0(\CC^\bullet (Z_1))$.
Suppose $z\in H^0(\CC^\bullet (Z_1))$. Then 
$z\in\mathrm{Im}(\psi_1)$ implies $\psi_1(z)=z$; and
$z\in\mathrm{Ker}(\psi_p)$ for all $p>1$ implies
$z\in H^0(\CC^\bullet (Z_0))$. Therefore,
$H^1(\CC^\bullet(Z))=0$. \qed
\end{example}

\subsection{}

Let $V$ be a $\A_{n,\la,\nu}$-module. 
For each $\jj\in I^n$, define 
$$Z_\jj(J):= V(\jj, \Delta(\jj)\setminus J) \quad
\mbox{ and } \quad \psi_{J,p} := \pi_{\jj,p}
\qquad \mbox{ for }
J\subset \Delta(\jj),\ p\in \Delta(\jj)\setminus J.$$
By Lemma \ref{firstlemma}(iv), $(Z_\jj,\psi)$ is 
a commutative $\Delta(\jj)$-cube (over the category of
$\B$-modules).
Define the complex of $\B\rtimes \k[S_n]$-modules
$$ \CC^\bullet(V) := \bigoplus_{\jj\in I^n}
\CC^\bullet (Z_\jj). $$
Thus,
$$ \CC^r (V) =
\bigoplus_{\jj\in I^n}
\bigoplus_{\substack{D\subset \Delta(\jj) \\
|D|=|\Delta(\jj)|-r}} V(\jj,D) \ot_\Z \det(\Delta(\jj)\setminus D), 
\qquad r=0, \ldots, n. $$
We remark that $S_n$ acts diagonally. Observe that 
\begin{equation} \label{h0}
F_i(V) = H^0(\CC^\bullet(V)). 
\end{equation}

We have the following results when $\nu=0$.
\begin{proposition} \label{highercohomology}
Let $\la\in B$ and assume $\la_i$ is invertible in $\k$.
Let $V$ be a $\A_{n,\la, 0}$-module.
Then $H^r (\CC^\bullet (V)) = 0$ for all $r>0$.
\end{proposition}
\begin{proof}
Fix $\jj\in I^n$. 
Let $Z'=V(\jj,\Delta(\jj))$. For each $p\in \Delta(\jj)$,
define an endomorphism $\psi_p$ of $Z'$ by
$\psi_p = \la_i^{-1} \mu_{\jj,p}\pi_{\jj,p}$.
Suppose $D=\Delta(\jj)\setminus \{q_1,\ldots,q_r\}$
and $p\in D$. Then
by Lemma \ref{firstlemma}, the following diagram commutes:
$$\xymatrix{
V(\jj, D) 
\ar[rrrr]^{\la_i^{-r}\mu_{\jj,q_1}\cdots\mu_{\jj,q_r}}
\ar[d]_{\pi_{\jj,p}} &&&& 
\mathrm{Im}(\psi_{q_1}\cdots\psi_{q_r})
\ar[d]^{\psi_p} \\
V(\jj, D\setminus\{p\})
\ar[rrrr]^{\la_i^{-(r+1)}\mu_{\jj,p}\mu_{\jj,q_1}\cdots\mu_{\jj,q_r}}
&&&& \mathrm{Im}(\psi_p\psi_{q_1}\cdots\psi_{q_r}) }$$
Moreover, the horizontal maps in the above diagram are 
isomorphisms.
Hence, the proposition is immediate 
from the claim in Example \ref{idempotent}.
\end{proof}

The Grothendieck group of an abelian category $\mathscr Z$
is an abelian group with generators $[Z]$, for all
objects $Z$ of $\mathscr Z$, and defining relations
$[Z] = [Z']+[Z'']$ for all short exact sequences
$0\to Z' \to Z \to Z'' \to 0$.

\begin{corollary} \label{eulerpoincare}
Let $\la\in B$ and assume $\la_i$ is invertible in $\k$.
Let $V$ be a $\A_{n,\la, 0}$-module.
Then in the Grothendieck group of the category of
$\B\rtimes \k[S_n]$-modules, we have
$$ F_i(V) = \sum_{r=0}^n (-1)^r
 \bigg[ \bigoplus_{\jj\in I^n}
  \bigoplus_{ \substack{D\subset\Delta(\jj) \\
                       |D|=|\Delta(\jj)|-r} }
  V(\jj, D) \ot_\Z \det(\Delta(\jj)\setminus D) \bigg] . $$
\end{corollary}
\begin{proof}
This follows from (\ref{h0}) and
Proposition \ref{highercohomology}
by the Euler-Poincar\'e principle.
\end{proof}

The author does not know what happens if $\nu\neq 0$;
but see Proposition \ref{flat2} below.

We conjecture that in general $H^r(\CC^\bullet(V))$
are the right derived functors of $F_i$.
One can also similarly define a complex using the maps
$\mu_{\jj,p}$ instead of $\pi_{\jj,p}$.

\subsection{}
In this subsection, we let $\k:=\C$.
We shall determine the set $\Lambda_i$.

First, we recall some standard results on the 
representation theory of symmetric groups;
see \cite[\S2.2]{EM} and \cite[\S2.4]{M1}.
For a Young diagram $\mu$ corresponding to a partition
of $n$, we write $\pi_\mu$ for the associated irreducible
representation of $S_n$.
For a cell $j$ in $\mu$, we let $\mathbf c(j)$ be the 
signed distance from $j$ to the diagonal.
The content $\mathbf c(\mu)$ of $\mu$ is the sum 
of $\mathbf c(j)$ over all cells $j$ in $\mu$.

Denote by $S_{n-1}$ the subgroup of $S_n$ which fixes $1$.
It is known that $\pi_\mu\big|_{S_{n-1}} = \bigoplus \pi_{\mu-j}$,
where the direct sum is taken over all corners $j$ of $\mu$,
and $\mu-j$ is the Young diagram obtained from $\mu$ be removing
the corner $j$.

\begin{lemma} \label{corner}
Let $C=s_{12}+s_{13}+\cdots+s_{1n}$.

{\rm (i)} The element $C$ acts on $\pi_{\mu-j}$ by the
scalar $\mathbf c(j)$, for each corner $j$ of 
the Young diagram $\mu$.

{\rm (ii)} The element $C$ acts as a scalar on $\pi_\mu$
if and only if $\mu$ is a rectangle. If the rectangle has
height $a$ and width $b$, then the scalar is $b-a$.
\end{lemma}

We omit the proof of the lemma, which can be found in
\cite[\S2.2]{EM} and \cite[\S2.4]{M1}.

%

\begin{proposition} \label{lambdai}
We have $$\Lambda_i = \{(\la,\nu)\in B\times \C
\mid \la_i \pm p\nu \neq 0 \mbox{ for } 
p = 0, 1, \ldots, n-1\}.$$
\end{proposition}
\begin{proof}
Let $r\in [1,n]$.
The element $\la_i + \nu\sum_{m=2}^r s_{1m}$ is invertible
in $\C[S_r]$ if and only if its eigenvalues on the irreducible
representations of $S_r$ are nonzero.
By Lemma \ref{corner}(i), the eigenvalues are
the numbers $\la_i + \nu\mathbf c(j)$, where $j$ is
a corner in a Young diagram.
The numbers $\mathbf c(j)$
which occur are $0,\pm 1, \ldots, \pm(n-1)$.
\end{proof}

\subsection{}
In this subsection, we let $\k:=\C[[U]]$, where $U$ is a finite
dimensional vector space over $\C$.
Let $\m$ be the unique maximal ideal of $\k$.
If $V$ is a $\k$-module, we write $\overline V$ for $V/\m V$.
A $\A_{n,\la,\nu}$-module $V$ is a flat formal deformation
of a $\overline{\A_{n,\la,\nu}}$-module
$V_0$ if $V\cong V_0[[U]]$ as $\k$-modules, and
there is a given isomorphism $\overline V \cong V_0$ of
$\overline{\A_{n,\la,\nu}}$-modules.

For any $\k$-module $V$, its $\m$-filtration is the
decreasing filtration
$V\supset \m V \supset \m^2 V \supset \ldots$.
We define
$$\mathrm{Gr}_\m V
:= \prod_{h=0}^{\infty} \frac{\m^h V}{\m^{h+1}V}.$$

Let us also introduce the following notations.
We shall write $\widetilde \k$ for $\C[U]$,
$\widetilde B$ for $\oplus_{i\in I} \widetilde \k$,
and $\widetilde E$ for the free $\widetilde \k$-module with
basis the set of edges $\{a\in \overline Q\}$.
Furthermore, let
$$ \widetilde\B := \widetilde B^{\otimes n},  \qquad
\widetilde \E := \bigoplus_{1\leq \ell \leq n}
\widetilde B^{\otimes (\ell-1)} \otimes \widetilde E \otimes
\widetilde B^{\otimes (n-\ell)}. $$
For any $\la\in\widetilde B$ and $\nu\in\C[U]$,
we write $\widetilde{\A_{n,\la,\nu}}$ for the $\widetilde \B$-algebra
defined as the quotient of $T_{\widetilde \B}\widetilde \E \rtimes
\widetilde \k[S_n]$ by the relations (i) and (ii) in 
Definition \ref{algebra}.

\begin{lemma} \label{pbw}
Let $Q$ be a connected quiver without edge-loop, such
that $Q$ is not a finite Dynkin quiver. Assume
$\la\in\widetilde B$ and $\nu\in \C[U]$.
Then $\mathrm{Gr}_\m \A_{n,\la,\nu} \cong
\overline{\A_{n,\la,\nu}}[[U]]$ as algebras over $\k$.
\end{lemma}
\begin{proof}
The algebra $\A_{n,\la,\nu}$ has an increasing filtration
defined by setting elements of $\B\rtimes \k[S_n]$ to
be of degree $0$, and elements of $\E$ to be of degree $1$.
Similarly, $\widetilde{\A_{n,\la,\nu}}$
and $\overline{\A_{n,\la,\nu}}$ have increasing filtrations.
                                                                            
Let $S'$ be a basis for $\overline \E$, and let
$S$ be a set of words in the elements of $S'$
such that $S$ is a basis for $\overline{\A_{n,0,0}}$
over $\overline\B\rtimes \C[S_n]$.
It was proved in \cite[Theorem 2.2.1]{GG} (see also
\cite[Remark 2.2.6]{GG}) 
that, for any $\la_0\in\overline B$
and $\nu_0\in \C$, the natural map
$\overline{\A_{n,0,0}}\to
\mathrm{gr}\overline{\A_{n,\la_0,\nu_0}}$
is an isomorphism of graded algebras.
Hence, $S$ is a basis for
$\overline{\A_{n,\la_0,\nu_0}}$ over $\overline\B\rtimes \C[S_n]$.
                                                                            
We have the natural epimorphism
$$ \overline{\A_{n,0,0}} \ot_\C \C[U]
= \widetilde{\A_{n,0,0}}
\too \mathrm{gr}\widetilde{\A_{n,\la,\nu}}.$$
Thus, $S$ spans $\widetilde{\A_{n,\la,\nu}}$
as a module over $\widetilde \B\rtimes\widetilde \k[S_n]$.
If there is a linear relation over
$\widetilde \B\rtimes\widetilde \k[S_n]$
among elements of $S$ in $\widetilde{\A_{n,\la,\nu}}$,
then by evaluation at some point of $U$,
we obtain a linear relation over
$\overline \B\rtimes \C[S_n]$, a contradiction.
Hence, $S$ is a basis for
$\widetilde{\A_{n,\la,\nu}}$ over $\widetilde \B\rtimes\widetilde \k
[S_n]$.
It follows that $S$ is a basis for
$\A_{n,\la,\nu} = \widetilde{\A_{n,\la,\nu}} \ot_{\C[U]}\k$
over $\B\rtimes\k[S_n]$.
Therefore, $\A_{n,\la,\nu}
\cong \overline{\A_{n,\la,\nu}}[[U]]$
as $\k$-modules, and
$\mathrm{Gr}_\m \A_{n,\la,\nu} \cong
\overline{\A_{n,\la,\nu}}[[U]]$
as algebras.
\end{proof}
                                                                          
\begin{proposition} \label{flat}
Let $Q$ and $\la,\nu$ be as in Lemma \ref{pbw}.
Let $i\in I$ and suppose $(\la,\nu)\in \Lambda_i$.
Let the $\A_{n,\la,\nu}$-module $V$ be a flat formal deformation
of a $\overline{\A_{n,\la,\nu}}$-module $V_0$.
Then $F_i(V)$ is a flat formal deformation of $F_i(V_0)$.
\end{proposition}
\begin{proof}
Observe that $\mathrm{Gr}_\m V = V_0[[U]]$ as
$\overline{\A_{n,\la,\nu}}[[U]]$-modules.
We have
$$\mathrm{Gr}_\m (F_i(V)) \subset F_i(\mathrm{Gr}_\m V).$$
By Theorem \ref{equivalence},
$$ \mathrm{Gr}_\m V =\mathrm{Gr}_\m (F_i F_iV)
\subset F_i(\mathrm{Gr}_\m(F_i V)) \subset
F_i(F_i(\mathrm{Gr}_\m V)) = \mathrm{Gr}_\m V .$$
Hence, $\mathrm{Gr}_\m(F_i(V)) = F_i(\mathrm{Gr}_\m V)
= (F_i V_0)[[U]]$ as $\overline{\A_{n,r_i\la,\nu}}[[U]]$-modules,
which implies
$F_i(V) \cong (F_iV_0)[[U]]$ as $\k$-modules,
and $\overline{F_i(V)}= F_i(V_0)$ as
$\overline{\A_{n,r_i\la,\nu}}$-modules.
\end{proof}

\subsection{}
In this subsection, we let $\k:=\C[U]$, the ring of regular
functions on a connected smooth affine variety $U$.
For any point $u\in U$, we denote by $\m_u$ the maximal 
ideal of functions vanishing at $u$, and if 
$V$ is a $\k$-module, then let $V^u:= V/\m_u V$. 
We shall write $\overline\B$ for $\B^u$.

\begin{proposition} \label{flat2}
Assume $Q$ is a connected quiver without edge-loop, such that
$Q$ is not a finite Dynkin quiver. Let $i\in I$ and
$(\la,\nu)\in \Lambda_i$.
Let $V$ be a $\A_{n,\la,\nu}$-module, finitely generated
over $\k$.
Suppose $V$ is a flat $\k$-module.
Then we have the following.

{\rm (i)}
$F_i(V)$ is a flat $\k$-module.

{\rm (ii)}
If $\nu$ vanishes at a point $o\in U$, then
for any point $u\in U$, we have
$$ F_i(V^u) = \sum_{r=0}^n (-1)^r
 \bigg[ \bigoplus_{\jj\in I^n}
  \bigoplus_{ \substack{D\subset\Delta(\jj) \\
                       |D|=|\Delta(\jj)|-r} }
  V^u(\jj, D) \ot_\Z \det(\Delta(\jj)\setminus D) \bigg] $$
in the Grothendieck group of the category of
$\overline\B \rtimes \C[S_n]$-modules.
\end{proposition}
\begin{proof}
(i) By Corollary \ref{fiber} and
Proposition \ref{flat}, $F_i(V)$ is locally flat 
at all the points of $U$. Hence, it is flat over $U$.

(ii) Since $V$ and $F_i(V)$ are flat over $U$, we have 
isomorphisms of $\overline\B\rtimes \C[S_n]$-modules:
$V^o \cong V^u$ and $F_i(V)^o \cong F_i(V)^u$.
Hence, the required formula follows from 
Corollary \ref{fiber} and Corollary \ref{eulerpoincare}.
\end{proof}

\section{\bf Symplectic reflection algebras for wreath products}

\subsection{}
Let $L$ be a 2-dimensional complex vector space,
and $\omega_L$ a nondegenerate symplectic form on $L$.
Let $\Gamma$ be a finite subgroup of $Sp(L)$, and
let $\GG:=S_n\ltimes\Gamma^n$.
Let $\mathscr L:=L^{\oplus n}$. 
For any $\ell\in [1,n]$ and $\g\in\Gamma$, we will write
$\g_\ell\in\GG$ for $\g$ placed in the $\ell$-th factor $\Gamma$.
Similarly, for any $u\in L$, we will write 
$u_\ell\in \mathscr L$ for $u$ placed in the $\ell$-th factor $L$.
Fix a symplectic basis $\{x,y\}$ for $L$.

Let $t,k\in\C$.
Denote by $Z\Gamma$ the center of the group algebra $\C[\Gamma]$.
Let $$c = \sum_{\g \in \Gamma\smallsetminus \{1\}} 
c_{\g}\cdot\g \in Z\Gamma, \quad\mbox{ where } c_\g\in\C.$$
The symplectic reflection algebra $\hh_{t,k,c}(\GG)$,
introduced in \cite{EG}, is the quotient of
$T\mathscr{L}\rtimes \C[\GG]$ by the following relations:
\begin{align*}
[x_\ell, y_\ell]
=& t\cdot 1+ \frac{k}{2} \sum_{m\neq \ell}\sum_{\g\in\Gamma}
s_{\ell m}\g_{\ell}\g_{m}^{-1} 
+ \sum_{\g\in\Gamma\smallsetminus\{1\}} c_{\g}\g_{\ell} , 
\qquad \forall  \ell\in[1,n]; \\
[u_\ell,v_m]=& -\frac{k}{2} \sum_{\g\in\Gamma} \omega_{L}(\g u,v)
s_{\ell m}\g_{\ell}\g_{m}^{-1},
\qquad \forall  u,v\in L,\ \ell,m\in [1,n],\ \ell\neq m. 
\end{align*}

Let $N_i$ be the irreducible representation of $\Gamma$
corresponding to the vertex $i\in I$ of $Q$ 
(where $Q$ is associated to $\Gamma$ by the McKay correspondence)
and let
$f_i\in \mathrm{End} N_i$ be a primitive idempotent.
We have $\C\Gamma=\bigoplus_{i\in I}  \mathrm{End} N_i$.
Let $f:=\sum_{i\in I}f_i\in\C\Gamma$.
The element $f^{\ot n}$ is an idempotent in $\C[\Gamma^n]
=(\C\Gamma)^{\ot n}$.
It was proved in \cite[Theorem 3.5.2]{GG} that the algebra
$f^{\ot n} \hh_{t,k,c}(\GG) f^{\ot n}$ is isomorphic
to the algebra $\A_{n,\la,\nu}$ for the quiver $Q$.
In particular, $\hh_{t,k,c}(\GG)$ is Morita equivalent
to $\A_{n,\la,\nu}$.
The parameter $\la_i$ is the trace of
$t\cdot 1+c$ on $N_i$, and the parameter
$\nu$ is $\frac{k|\Gamma|}{2}$.
We shall reformulate and prove the main results of
\cite{EM}, \cite{M1} and \cite{M2} in terms of the algebra
$\A_{n,\la,\nu}$ via this Morita equivalence.
We believe the results are more transparent in our
reformulation.

\subsection{}
In this subsection, we let $\k:=\C[[U]]$, where $U$ is a finite
dimensional vector space over $\C$. 
Let $\m$ be the unique maximal ideal of $\k$. Recall that
if $V$ is a $\k$-module, we write $\overline V$ for $V/\m V$.

If $V=\bigoplus_{i\in I}V_i$ is an $I$-graded complex vector space,
then its dimension vector is the element
$(\dim V_i)_{i\in I}\in\Z^I$.
Let $\N_i$ be the complex vector space with dimension vector
$\epsilon_i$.

Let $\vec n=(n_1, \ldots, n_r)$ be a partition of $n$.
Let $X=X_1 \otot X_r$ be a simple module of
$S_{\vec n}:=S_{n_1}\times \cdots \times S_{n_r}
\subset S_n$.
Let $\{i_1, \ldots, i_r\}$ be a set of $r$ \emph{distinct}
vertices of $Q$, and let 
$\N=\N_{i_1}^{\ot n_1} \otot \N_{i_r}^{\ot n_r}$.
Then $X\ot\N$ is a simple module of
$\overline \B \rtimes \C[S_{\vec n}]$.
We write $X\ot\N\uparrow$ for the induced module
$\mathrm{Ind}_{\overline{\B} \rtimes \C[S_{\vec n}]}
^{\overline{\B} \rtimes \C[S_n]} (X\ot \N)$
of $\overline{\B} \rtimes \C[S_n]$.
It is known that
any simple $\overline{\B} \rtimes \C[S_n]$-module
is of the  form $X\ot \N\uparrow$ 
(see \cite{Mac}: paragraph after (A.5)). 
We have
$$ X\ot\N\uparrow = \bigoplus_\sigma \sigma(X\ot\N),$$
where $\sigma$ runs over a set of left coset representatives
of $S_{\vec n}$ in $S_n$.
                                                                             
The following lemma is equivalent to \cite[Theorem 4.1]{M1}.
                                                                             
\begin{lemma} \label{trivial}
Assume $Q$ has no edge-loop.
Let the $\A_{n,\la,\nu}$-module $V$ be a flat formal deformation
of the $\overline{\A_{n,\la,\nu}}$-module $\overline V$.
If $\overline V$ is simple as a $\overline\B\rtimes\C[S_n]$-module,
then all elements of $\E$ must act by $0$ on $V$.
\end{lemma}
\begin{proof}
Since the algebra $\overline\B\rtimes\C[S_n]$ is semisimple,
$V$ must be of the form $(X\ot\N\uparrow)[[U]]$
(as $\B\rtimes\k[S_n]$-modules).
Let $(j_1,\ldots, j_n)\in I^n$ and $\sigma\in S_n$.
If $j_m = j_{\sigma(m)}$ for all $m\neq \ell$, then
$j_\ell = j_{\sigma(\ell)}$.
It follows that since $Q$ has no edge-loop, 
$\E_\ell$ must act by $0$ on $\N$, hence $\E$ acts
by $0$ on $V$.
\end{proof}

The next result is equivalent to \cite[Theorem 3.1]{M1}.

\begin{theorem} \label{extend}
Assume $Q$ has no edge-loop and $\nu\neq 0$.
The $\B\rtimes \k[S_n]$-module $(X\ot\N\uparrow)[[U]]$
extends to a $\A_{n,\la,\nu}$-module if and only if
the following conditions are satisfied:
\begin{itemize}
\item[{\rm (i)}]
For all $\ell\in[1,r]$, the simple module $X_\ell$
of $S_{n_\ell}$ has rectangular Young diagram, of
size $a_\ell \times b_\ell$.
\item[{\rm (ii)}]
No two vertices in $\{i_1, \ldots, i_r\}$ are joined
by an edge in $Q$.
\item[{\rm (iii)}]
For all $\ell\in[1,r]$, one has
$\la_{i_\ell} = (a_\ell-b_\ell)\nu$.
\end{itemize}
\end{theorem}
\begin{proof}
Suppose $(X\ot\N\uparrow)[[U]]$
extends to a $\A_{n,\la,\nu}$-module.
By Lemma \ref{trivial}, the elements of $\E$ must act by $0$.
Hence, by Lemma \ref{corner}(ii) and
the relations of type (i) in Definition \ref{algebra},
the Young diagram of each $X_\ell$ must be a rectangle,
of size $a_\ell\times b_\ell$ say, and 
$\la_{i_\ell} = (a_\ell-b_\ell)\nu$.
By the relations of type (ii) in Definition \ref{algebra},
no two vertices in $\{i_1, \ldots, i_r\}$ can be joined
by an edge in $Q$.

Conversely, suppose the conditions are satisfied.
Then it is clear (using Lemma \ref{corner}(ii) again) 
that if we let the elements of $\E$ act by $0$, the
relations in Definition \ref{algebra} hold.
\end{proof}

From now on, we assume that $Q$ is an affine Dynkin
quiver of type ADE, but not of type $\mathrm{A}_0$.
Let $\delta=(\delta_i)_{i\in I}\in \Z^I$ be the 
minimal positive imaginary root of $Q$.
We have $\delta_i = \dim N_i$.

Let $\la_0\in\overline B$, and assume $\la_0\cdot\delta\neq 0$. 
We shall write
$\Pi_{\la_0}$ for $\overline{\Pi_{\la_0}}$,
and $\A_{n,\la_0,0}$ for $\overline{\A_{n,\la_0,0}}$.

Let $\Sigma_{\la_0}$ be the set of dimension vectors
of finite dimensional simple $\Pi_{\la_0}$-modules.
By \cite[Lemma 7.2]{CBH} and \cite[Theorem 7.4]{CBH},
there exists an element $\la^+\in\overline B$ and an
element $w\in W$ such that:
\begin{itemize}
\item
$w$ is an element of minimal length with $w(\la_0)=\la^+$;
\item
$w\Sigma_{\la_0} = \{\epsilon_i \mid \la^{+}_i=0\}$.
\end{itemize}
By the minimality of its length, we can write
$w=s_{j_h}\cdots s_{j_1}$ for some
$j_1, \ldots, j_h\in I$, such that
$\big(r_{j_g} \cdots r_{j_1}(\la_0)\big)_{j_g}\neq 0$
for all $g\in [1,n]$.
Let $F_w$ be the composition of functors
$F_{j_h}\cdots F_{j_1}$, and
$F_{w^{-1}}$ be the composition of functors
$F_{j_1}\cdots F_{j_h}$.

Let $Y_1, \ldots, Y_r$ be a collection of pairwise non-isomorphic
finite dimensional simple modules of $\Pi_{\la_0}$, and let
$Y=Y_1^{\ot n_1}\otot Y_r^{\ot n_r}$.
Then $X\ot Y$ is a simple module of
$\Pi_{\la_0}^{\ot n} \rtimes \C[S_{\vec n}]$.
We write $X\ot Y\uparrow$ for the induced module
$\mathrm{Ind}_{\Pi_{\la_0}^{\ot n} \rtimes \C[S_{\vec n}]}
^{\Pi_{\la_0}^{\ot n} \rtimes \C[S_n]} (X\ot Y)$
of $\A_{n,\la_0,0}$.
By \cite{Mac} (paragraph after (A.5)),
it is known that any finite dimensional simple 
$\A_{n,\la_0,0}$-module is of the form $X\ot Y\uparrow$.

The following theorem and its proof was explained 
to the author by Pavel Etingof.
\begin{theorem} \label{deform}
Let $\la\in B$.
Assume that $\la_i\in U$ for all $i\in I$, 
and $0\neq \nu\in U$.
The $\A_{n,\la_0,0}$-module 
$X\ot Y\uparrow$ has a flat formal deformation
to a $\A_{n,\la_0+\la,\nu}$-module if and only if 
the following conditions are satisfied:
\begin{itemize}
\item[{\rm (i)}]
For all $\ell\in[1,r]$, the simple module $X_\ell$
of $S_{n_\ell}$ has rectangular Young diagram, of
size $a_\ell \times b_\ell$.
\item[{\rm (ii)}]
We have $\mathrm{Ext}^1_{\Pi_{\la_0}}(Y_\ell,Y_m)=0$ 
for any $\ell\neq m$.
\item[{\rm (iii)}]
For all $\ell\in[1,r]$, one has
$\la \cdot \al_\ell = (a_\ell-b_\ell)\nu$, where
$\alpha_\ell$ is the dimension vector of $Y_\ell$.
\end{itemize}
The deformation is unique when it exists.
\end{theorem}
\begin{proof}
Let $\la^+$, $w$, $F_w$, and $F_{w^{-1}}$ be as 
defined above.
We claim that
$\big(r_{j_g}\cdots r_{j_1} (\la_0 +\la), \nu\big)
\in \Lambda_{j_g}$ for $g=1,\ldots,h$. 
To see this, it is enough to show that 
$\big(r_{j_g}\cdots r_{j_1} (\la_0+\la)\big)_{j_g} +
 \nu C $ has an inverse in $\k[S_N]$, for any given $C\in \k[S_N]$.
This is equivalent to solving a system
of $N!$ linear equations in $N!$ variables,
whose associated matrix is of the form
$\big(r_{j_g}\cdots r_{j_1} (\la_0+\la)\big)_{j_g} Id_{N!} +
 \nu M$ for some matrix $M$.
Since $\la$ and $\nu$ are $0$ modulo $\m$, the
determinant of this matrix is nonzero modulo $\m$, and so
it is invertible in $\k$. Hence, the matrix is
invertible. This proved  our claim.

Now define $i_\ell\in I$ by $\epsilon_{i_\ell} = w(\al_\ell)$.
We have $\la_0\cdot \al_\ell = \la^+ \cdot \epsilon_{i_\ell} =0$.
By \cite[Theorem 5.1]{CBH},
we have $F_w(X\ot Y\uparrow)=X\ot\N\uparrow$.

Suppose the conditions in the theorem are satisfied.
Then by Theorem \ref{extend}, the $\B\rtimes \k[S_n]$-module
$M := (X\ot\N\uparrow)[[U]]$ is a $\A_{n,\la^+ +w(\la),\nu}$-module 
(where elements of $\E$ act by $0$).
Hence, by Proposition \ref{flat}, the $\A_{n,\la_0+\la,\nu}$-module
$F_{w^{-1}} (M)$ is a flat formal deformation of $X\ot Y\uparrow$.

Conversely, suppose a $\A_{n,\la_0+\la,\nu}$-module $V$ is a
flat formal deformation of $X\ot Y\uparrow$. 
Then by Proposition \ref{flat}, the 
$\A_{n,\la^+ +w(\la),\nu}$-module $F_w(V)$ is a 
flat formal deformation 
of $X\ot\N\uparrow$.
We have $F_w(V) = (X\ot\N\uparrow)[[U]]$
as $\B\rtimes\k[S_n]$-modules.
It follows from Theorem \ref{extend} that the conditions 
in the theorem must hold.
Moreover, by Lemma \ref{trivial}, the elements of $\E$ must
act by $0$ on $F_w(V)$, so $F_w(V)$ is the unique flat 
formal deformation
of $X\ot\N\uparrow$. This implies that $V$ is the unique 
flat formal deformation of $X\ot Y\uparrow$.
\end{proof}

The sufficiency of the conditions in Theorem \ref{deform}
was first proved in \cite[Theorem 1.3(i)]{M2};
in the special case where the partition is $\vec n=(n)$,
it was first proved in \cite[Theorem 3.1(i)]{EM}.

\subsection{}
Let $\la_0$, $\la^+$, $w$, $F_{w^{-1}}$, 
$j_1, \ldots, j_h$, and $X\ot Y\uparrow$
be as defined in the previous subsection.
In particular, $\la_0\cdot\delta\neq 0$.
Assume that conditions (i) and (ii) of Theorem \ref{deform} hold.

Let $U$ be a finite dimensional complex vector space.
Let $\la$ and $\nu$ be regular functions on $U$ such that
condition (iii) of Theorem \ref{deform} hold. Moreover,
assume there is a point $o\in U$ such that
$\la$ specializes to $\la_0$ at $o$, and
$\nu$ vanishes at $o$.

Define $i_\ell\in I$ by $\epsilon_{i_\ell} = w(\al_\ell)$,
and let $X\ot\N\uparrow$ be as defined in the previous subsection.
We have $\la^+_{i_\ell}=\la_0\cdot\alpha_\ell =0$.

Let $U'$ be the Zariski open subset of $U$ defined by 
$\big(r_{j_g} \cdots r_{j_1}(\la)\big)_{j_g}
\pm p\nu \neq 0$
for all $g\in [1,n]$ and $p=0, \ldots, n-1$.
Since $o\in U'$, the set $U'$ is nonempty.

Let $\k:=\C[U']$ be the ring of regular functions on $U'$. 
For any point $u\in U'$, let $\m_u$
denote the maximal ideal of $\k$ consisting of functions
vanishing at $u$.
If $V$ is a $\k$-module, we write $V^u$ for $V/\m_u V$.
We write $\overline\B$ for $\B^u$.

The proof of the following theorem  
is similar to the proof of Theorem \ref{deform}.

\begin{theorem} \label{deform2}
There exists a $\A_{n,\la,\nu}$-module $V$ such that:
\begin{itemize}
\item[{\rm (i)}]
$V^o=X\ot Y\uparrow$ as $\A_{n,\la_0,0}$-modules, 
and $V$ is flat over $U'$.
\item[{\rm (ii)}]
For any point $u\in U'$,
$V^u$ is a finite dimension simple $\A_{n,\la,\nu}^u$-module,
isomorphic to $X\ot Y\uparrow$ as a 
$\overline\B \rtimes \C[S_n]$-module.
\end{itemize}
\end{theorem}
\begin{proof}
Let the elements of $\E$ act by $0$ on the
$\B\rtimes\k[S_n]$-module $(X\ot \N \uparrow)\ot_\C \k$.
It follows from Lemma \ref{corner}(ii) that
$(X\ot \N \uparrow)\ot_\C \k$
is a $\A_{n,w(\la),\nu}$-module.

By Proposition \ref{lambdai}, we have
$\big( r_{j_g} \cdots r_{j_1}(\la) ,\nu \big) \in\Lambda_{j_g}$
for $g=1,\ldots,h$.
Let $V$ be the $\A_{n,\la,\nu}$-module
$F_{w^{-1}}\big((X\ot \N \uparrow)\ot_\C \k\big)$.
By Corollary \ref{fiber}, we have $V^o=X\ot Y \uparrow$.
Moreover, by Theorem \ref{equivalence},
$V^u$ is a simple $\A_{n,\la,\nu}^u$-module for any $u\in U'$.
By Proposition \ref{flat2}(i), $V$ is flat over $U'$.
Hence, $V^u$ is isomorphic to $V^o$
as $\overline\B\rtimes\C[S_n]$-modules.
\end{proof}
We remark that the set $U'$ was not specified precisely in
\cite[Theorem 3.1(iii)]{EM} and \cite[Theorem 1.3(iii)]{M2}.

Let us also mention that there may exists finite dimensional simple 
modules of $\hh_{t,k,c}(\GG)$ (for complex parameters) which 
cannot be deformed to a flat family as $k$ varies;
see \cite[\S4]{Ch} where this happens.

The assumption that $\la_0\cdot\delta\neq 0$ is equivalent
to the condition that $t\neq 0$ for the symplectic
reflection algebra $\hh_{t,k,c}(\GG)$.
When $t=0$, the representation theory of the symplectic
reflection algebra is remarkably different;
see \cite{CBH}, \cite{EG}, and \cite{GS}.

\section*{\bf Acknowledgments}

I am very grateful to Pavel Etingof for 
patiently explaining Theorem \ref{deform} and its proof to me,
and for other useful discussions.
I also thank Iain Gordon for his comments.
This work was partially supported by NSF grant DMS-0401509.

{\footnotesize {

}}


\begin{thebibliography}{APK}

\bibitem[BEG1]{BEG1} Y. Berest, P. Etingof, V. Ginzburg,
{\it Cherednik algebras and differential operators on 
quasi-invariants},
Duke Math. J. {\bf 118}  (2003), no. 2, 279--337,
{\tt math.QA/0111005}.

\bibitem[BEG2]{BEG2} Y. Berest, P. Etingof, V. Ginzburg,
{\it Morita Equivalence of Cherednik Algebras},
J. Reine Angew. Math. {\bf  568} (2004), 81--98,
{\tt math.QA/0207295}.

\bibitem[BGP]{BGP} I.N. Bernstein, I.M. Gelfand, V.A. Ponomarev,
{\it Coxeter functors, and Gabriel's theorem}, 
Russian Math. Surveys {\bf 28} (1973), no. 2, 17--32.

\bibitem[CBH]{CBH} W. Crawley-Boevey, M.P. Holland,
{\it Noncommutative deformations of Kleinian singularities},
Duke Math. J. {\bf 92} (1998), no. 3, 605--635.

\bibitem[Ch]{Ch} T. Chmutova,
{\it Representations of the rational Cherednik algebras of 
dihedral type}, preprint,
{\tt math.RT/0405383}.

\bibitem[EG]{EG} P. Etingof, V. Ginzburg,
{\it Symplectic reflection algebras, Calogero-Moser space, and 
deformed Harish-Chandra homomorphism},
Invent. Math. {\bf 147} (2002), 243-348, {\tt{math.AG/0011114}}.

\bibitem[EM]{EM} P. Etingof, S. Montarani,
{\it Finite dimensional representations of symplectic reflection
algebras associated to wreath products},
preprint, {\tt{math.RT/0403250}}.

\bibitem[GG]{GG} W.L. Gan, V. Ginzburg,
{\it Deformed preprojective algebras and symplectic reflection 
algebras for wreath products},
J. Algebra {\bf 283} (2005), no. 1, 350--363,
{\tt{math.QA/0401038}}.

\bibitem[GS]{GS} I. Gordon, S.P. Smith,
{\it Representations of symplectic reflection algebras and
resolutions of deformations of symplectic quotient singularities},
Math. Ann. {\bf 330} (2004), no. 1, 185--200, 
{\tt{math.RT/0310187}}.

\bibitem[Kh]{Kh} M. Khovanov,
{\it A categorification of the Jones polynomial},
Duke Math. J. {\bf 101} (2000), no. 3, 359--426,
{\tt{math.QA/9908171}}.

\bibitem[Mac]{Mac}  I.G. Macdonald,
{\it  Polynomial functors and wreath products},
J. Pure Appl. Algebra 18 (1980), no. 2, 173--204.

\bibitem[M1]{M1} S. Montarani,
{\it On some finite dimensional representations of symplectic 
reflection algebras associated to wreath products},
preprint, {\tt{math.RT/0411286}}.

\bibitem[M2]{M2} S. Montarani,
{\it Finite dimensional representations of symplectic reflection 
algebras associated to wreath products II},
preprint, {\tt{math.RT/0501156}}.

\bibitem[Na]{Na} H. Nakajima,
{\it Reflection functors for quiver varieties 
and Weyl group actions},
Math. Ann. {\bf 327} (2003), no. 4, 671--721. 

\bibitem[We]{We} C. Weibel,
{\it An introduction to homological algebra}, 
Cambridge Studies in Advanced Mathematics {\bf 38}, 
Cambridge University Press, Cambridge, 1994. 

\end{thebibliography}
\end{document}